\documentclass[final,3p]{elsarticle}
 \usepackage{graphics}
 \usepackage{graphicx}
 \usepackage{epsfig}
\usepackage{amssymb}
 \usepackage{amsthm}
 \usepackage{lineno}
 \usepackage{amsmath}
   \numberwithin{equation}{section}
\usepackage{mathrsfs}

\newtheorem{thm}{Theorem}[section]
\newtheorem{cor}[thm]{Corollary}
\newtheorem{lem}[thm]{Lemma}

\newtheorem{defn}[thm]{Definition}

\biboptions{sort&compress,square}
\journal{some mathematical journal}
\begin{document}
\begin{frontmatter}
\author[rvt1]{Sining Wei}
\ead{weisn835@nenu.edu.cn}
\author[rvt2]{Yong Wang\corref{cor2}}
\ead{wangy581@nenu.edu.cn}

\cortext[cor2]{Corresponding author.}

\address[rvt1]{School of Data Science and Artificial Intelligence, Dongbei University of Finance and Economics, \\
Dalian, 116025, P.R.China}
\address[rvt2]{School of Mathematics and Statistics, Northeast Normal University, Changchun, 130024, P.R.China}

\title{Perturbations of Dirac operators, spectral Einstein functionals and the Noncommutative residue}
\begin{abstract}
In this paper, we introduce the spectral Einstein functional for perturbations of Dirac operators on manifolds with boundary. Furthermore, we provide the proof of the Dabrowski-Sitarz-Zalecki type theorems associated with the spectral Einstein functionals for perturbations of Dirac operators, particularly in the cases of on 4-dimensional manifolds with boundary.
\end{abstract}

\end{frontmatter}
\section{Introduction}
\label{1}
An eminent spectral scheme, yielding geometric objects on manifolds, including residue, scalar curvature, and various scalar combinations of curvature tensors,
is the small-time asymptotic expansion of the (localised) trace of heat kernel \cite{PBG,FGV}. The theory exhibits profound and rich structures, spanning both physics and mathematics.
In the recent paper \cite{DL}, Dabrowski etc. introduced bilinear functionals of vector fields and differential forms,
the densities of which yield the  metric and Einstein spectral functionals on even-dimensional Riemannian manifolds,
and they obtained certain
values or residues of the (localised) zeta function of the Laplacian  arising from
the Mellin transform and the coefficients of this expansion.

Consider $E$ as a finite-dimensional complex vector bundle over a closed compact manifold $M$ with dimension $n$, the noncommutative residue of a pseudo-differential operator
$P\in\Psi DO(E)$ can be defined as
 \begin{equation}\label{1}
res(P):=(2\pi)^{-n}\int_{S^{*}M}{\rm trace}(\sigma_{-n}^{P}(x,\xi))\mathrm{d}x \mathrm{d}\xi,
\end{equation}
where $S^{*}M\subset T^{*}M$ represents the co-sphere bundle on $M$ and
$\sigma_{-n}^{P}$ is the component of order $-n$ of the complete symbol
 \begin{equation}\label{2}
\sigma^{P}:=\sum_{i}\sigma_{i}^{P}
\end{equation}
of $P$,  as defined in \cite{Ac,Wo,Wo1,Gu},
and the linear functional $res: \Psi DO(E)\rightarrow \mathbb{C }$
is, in fact, the unique trace (up to multiplication
by constants) on the algebra of pseudo-differential operators $\Psi DO(E)$.

In \cite{Co1}, Connes used the noncommutative residue to derive a conformal 4-dimensional Polyakov action analogy.
Connes  proved that the noncommutative residue on a compact manifold $M$ coincided with Dixmier's trace on pseudodifferential
operators of order -dim$M$\cite{Co2}. Furthermore, Connes claimed the noncommutative residue of the square of the inverse of the Dirac operator was proportioned to the Einstein-Hilbert action in \cite{Co2}.
Kastler provided a brute-force proof of this theorem in \cite{Ka}, while Kalau and Walze proved it in the normal coordinates system simultaneously in \cite{KW},
which is called the Kastler-Kalau-Walze theorem now.
Building upon the theory of the noncommutative reside  introduced by Wodzicki, Fedosov etc. \cite{FGLS} constructed a noncommutative
residue on the algebra of classical elements in Boutet de Monvel's calculus on a compact manifold with boundary of dimension $n>2$.
With elliptic pseudodifferential operators and  noncommutative
residue, it is natural way for investigating the Kastler-Kalau-Walze type theorem and
operator-theoretic explanation of the gravitational action on manifolds with boundary. Concerning Dirac operators and signature operators, Wang performed computations of the noncommutative residue and successfully demonstrated the Kastler-Kalau-Walze type theorem for manifolds with boundaries
\cite{Wa1,Wa3,Wa4}.

Jean-Michel Bismut \cite{JMB} previously established a local index theorem for Dirac operators on a Riemannian manifold $M$, specifically those associated with connections on $TM$ with non zero torsion.
In \cite{AT}, Ackermann and Tolksdorf demonstrated a generalized version of the well-known Lichnerowicz formula. This extended formula applies to for the square of the
most general Dirac operator with torsion $D_{T}$ on an even-dimensional spin manifold, associated to a metric connection with torsion.
In \cite{PS}, Pf$\ddot{a}$ffle and Stephan focused on compact Riemannian spin manifolds without boundary equipped with orthogonal connections,
and delved into the induced Dirac operators. Furthermore, Pf$\ddot{a}$ffle and Stephan investigated
orthogonal connections with arbitrary torsion on compact Riemannian manifolds,
they computed the spectral action for the induced
Dirac operators, twisted Dirac operators and Dirac operators of Chamseddine-Connes type, as detailed in their work\cite{PS1}. In \cite{WWY},  J. Wang, Y. Wang and C. Yang calculated the lower
dimensional  volume $\widetilde{{\rm Wres}}[\pi^+(D^{*}_{T})^{-p_{1}}\circ\pi^+D_{T}^{-p_{2}}]$ and derived a Kastler-Kalau-Walze
type theorems associated with Dirac operators with torsion on compact manifolds with  boundary. In \cite{JYT},
the authors generalized the results in \cite{DL}, \cite{Wa3}, \cite{PS} and obtained spectral functionals
 associated with Dirac operators with torsion on compact manifolds with  boundary. In \cite{AAA}, Wang considered the arbitrary perturbations of Dirac operators, and established the associated Kastler-Kalau-Walze theorem.

$\mathbf{The~aim ~of~ this~ paper}$ is to extend and generalize the results in \cite{DL}, \cite{JYT} and focusing on obtaining the spectral Einstein functional for perturbations of Dirac operators on compact manifolds with  boundary.  Specifically, for lower dimensional compact Riemannian manifolds
  with  boundary, we compute
the noncommutative residue $\widetilde{{\rm Wres}}\bigg[\pi^{+}\nabla^{\Psi}_{X}\nabla^{\Psi}_{Y}D_{\Psi}^{-2}
    \circ\pi^{+}D_{\Psi}^{-2}\bigg]$ and $\widetilde{{\rm Wres}}\bigg[\pi^{+}\nabla^{\Psi}_{X}\nabla_{Y}^{\Psi}D_{\Psi}^{-1}
    \circ\pi^{+}D_{\Psi}^{-3}\bigg]$ on 4 dimensional manifolds. Our main theorems are as follows.

\begin{thm}\label{thm1.1}
Let $M$ be an $4$-dimensional oriented
compact spin manifold with boundary $\partial M$, then we get the following equality:
\begin{eqnarray}
\label{1111}
&&\widetilde{{\rm Wres}}\bigg[\pi^{+}\nabla^{\Psi}_{X}\nabla^{\Psi}_{Y}D_{\Psi}^{-2}
    \circ\pi^{+}D_{\Psi}^{-2}\bigg]=
    \frac{4\pi^2}{3}\int_{M}EG(X,Y)vol_{g}
+\frac{\pi^2}{2}\int_{M}{\rm trace}\Big[c(X)\nabla^{S(TM)}_Y\big(c(\Psi)\big)\nonumber\\
&+&\nabla^{S(TM)}_Y\big(c(\Psi)\big) c(X)-c(Y)\nabla^{S(TM)}_X(c(\Psi))-\nabla^{S(TM)}_X(c(\Psi)) c(Y)\Big]dvol_{M}+\frac{1}{2}\int_{M}{\rm trace}[\frac{1}{4}s\nonumber\\
&+&\sum\limits^{n}_{j=1}\frac{1}{2}c(\Psi)c(e_j)
c(\Psi)c(e_j)-\big(c(\Psi)\big)^2]g(X,Y)dvol_{M}
+\int_{\partial M}\widetilde{\Phi},
\end{eqnarray}
where $\widetilde{\Phi}$ are defined in (\ref{795}).
\end{thm}
\begin{thm}\label{thm1.2}
Let $M$ be an $4$-dimensional oriented
compact spin manifold with boundary $\partial M$, then we get the following equality:
\begin{eqnarray}\label{1}
\label{1b2773}
&&\widetilde{{\rm Wres}}\bigg[\pi^{+}\nabla^{\Psi}_{X}\nabla_{Y}^{\Psi}D_{\Psi}^{-1}
    \circ\pi^{+}D_{\Psi}^{-3}\bigg]=\frac{4\pi^2}{3}\int_{M}EG(X,Y)vol_{g}
+\frac{\pi^2}{2}\int_{M}{\rm trace}\Big[c(X)\nabla^{S(TM)}_Y\big(c(\Psi)\big)\nonumber\\
&+&\nabla^{S(TM)}_Y\big(c(\Psi)\big) c(X)-c(Y)\nabla^{S(TM)}_X(c(\Psi))-\nabla^{S(TM)}_X(c(\Psi)) c(Y)\Big]dvol_{M}+\frac{1}{2}\int_{M}{\rm trace}[\frac{1}{4}s\nonumber\\
&+&\sum\limits^{n}_{j=1}\frac{1}{2}c(\Psi)c(e_j)
c(\Psi)c(e_j)-\big(c(\Psi)\big)^2]g(X,Y)dvol_{M}
+\int_{\partial M}\widehat{\Phi},
\end{eqnarray}
where $\widehat{\Phi}$ are defined in (\ref{1795}).
\end{thm}

The paper is organized in the following way. In Section \ref{section:2}, we review
some basic formulas related to the spectral Einstein functional for perturbations of Dirac operators
.
In Section \ref{section:3}, we prove the Dabrowski-Sitarz-Zalecki type theorem for $\widetilde{{\rm Wres}}\bigg[\pi^{+}\nabla^{\Psi}_{X}\nabla^{\Psi}_{Y}D_{\Psi}^{-2}
\circ\pi^{+}D_{\Psi}^{-2}\bigg]$ on 4-dimensional manifolds with boundary. In Section \ref{section:4},  we prove the Dabrowski-Sitarz-Zalecki type theorem $\widetilde{{\rm Wres}}\bigg[\pi^{+}\nabla^{\Psi}_{X}\nabla_{Y}^{\Psi}D_{\Psi}^{-1}
    \circ\pi^{+}D_{\Psi}^{-3}\bigg]$ on 4-dimensional manifolds with boundary.

\section{Boutet de
Monvel's calculus and the noncommutative residue}
\label{section:2}
 In this section, we recall some basic facts and formulas about Boutet de
Monvel's calculus and the definition of the noncommutative residue for manifolds with boundary, which will be used in the following. For more details, see Section 2 in \cite{Wa3}.\\
 \indent Let $M$ be a 4-dimensional compact oriented manifold with boundary $\partial M$.
We assume that the metric $g^{TM}$ on $M$ has the following form near the boundary,
\begin{equation}
\label{b1}
g^{M}=\frac{1}{h(x_{n})}g^{\partial M}+dx _{n}^{2},
\end{equation}
where $g^{\partial M}$ is the metric on $\partial M$ and $h(x_n)\in C^{\infty}([0, 1)):=\{\widehat{h}|_{[0,1)}|\widehat{h}\in C^{\infty}((-\varepsilon,1))\}$ for
some $\varepsilon>0$ and $h(x_n)$ satisfies $h(x_n)>0$, $h(0)=1$ where $x_n$ denotes the normal directional coordinate. Let $U\subset M$ be a collar neighborhood of $\partial M$ which is diffeomorphic with $\partial M\times [0,1)$. By the definition of $h(x_n)\in C^{\infty}([0,1))$
and $h(x_n)>0$, there exists $\widehat{h}\in C^{\infty}((-\varepsilon,1))$ such that $\widehat{h}|_{[0,1)}=h$ and $\widehat{h}>0$ for some
sufficiently small $\varepsilon>0$. Then there exists a metric $g'$ on $\widetilde{M}=M\bigcup_{\partial M}\partial M\times
(-\varepsilon,0]$ which has the form on $U\bigcup_{\partial M}\partial M\times (-\varepsilon,0 ]$
\begin{equation}
\label{b2}
g'=\frac{1}{\widehat{h}(x_{n})}g^{\partial M}+dx _{n}^{2} ,
\end{equation}
such that $g'|_{M}=g$. We fix a metric $g'$ on the $\widetilde{M}$ such that $g'|_{M}=g$.

Let the Fourier transformation $F'$  be
\begin{equation*}
F':L^2({\bf R}_t)\rightarrow L^2({\bf R}_v);~F'(u)(v)=\int_\mathbb{R} e^{-ivt}u(t)dt
\end{equation*}
and let
\begin{equation*}
r^{+}:C^\infty ({\bf R})\rightarrow C^\infty (\widetilde{{\bf R}^+});~ f\rightarrow f|\widetilde{{\bf R}^+};~
\widetilde{{\bf R}^+}=\{x\geq0;x\in {\bf R}\}.
\end{equation*}
\indent We define $H^+=F'(\Phi(\widetilde{{\bf R}^+}));~ H^-_0=F'(\Phi(\widetilde{{\bf R}^-}))$ which satisfies
$H^+\bot H^-_0$, where $\Phi(\widetilde{{\bf R}^+}) =r^+\Phi({\bf R})$, $\Phi(\widetilde{{\bf R}^-}) =r^-\Phi({\bf R})$ and $\Phi({\bf R})$
denotes the Schwartz space. We have the following
 property: $h\in H^+~$ (respectively $H^-_0$) if and only if $h\in C^\infty({\bf R})$ which has an analytic extension to the lower (respectively upper) complex
half-plane $\{{\rm Im}\xi<0\}$ (respectively $\{{\rm Im}\xi>0\})$ such that for all nonnegative integer $l$,
 \begin{equation*}
\frac{d^{l}h}{d\xi^l}(\xi)\sim\sum^{\infty}_{k=1}\frac{d^l}{d\xi^l}(\frac{c_k}{\xi^k}),
\end{equation*}
as $|\xi|\rightarrow +\infty,{\rm Im}\xi\leq0$ (respectively ${\rm Im}\xi\geq0)$ and where $c_k\in\mathbb{C}$ are some constants.\\
 \indent Let $H'$ be the space of all polynomials and $H^-=H^-_0\bigoplus H';~H=H^+\bigoplus H^-.$ Denote by $\pi^+$ (respectively $\pi^-$) the
 projection on $H^+$ (respectively $H^-$). Let $\widetilde H=\{$rational functions having no poles on the real axis$\}$. Then on $\tilde{H}$,
 \begin{equation}
 \label{b3}
\pi^+h(\xi_0)=\frac{1}{2\pi i}\lim_{u\rightarrow 0^{-}}\int_{\Gamma^+}\frac{h(\xi)}{\xi_0+iu-\xi}d\xi,
\end{equation}
where $\Gamma^+$ is a Jordan closed curve
included ${\rm Im}(\xi)>0$ surrounding all the singularities of $h$ in the upper half-plane and
$\xi_0\in {\bf R}$. In our computations, we only compute $\pi^+h$ for $h$ in $\widetilde{H}$. Similarly, define $\pi'$ on $\tilde{H}$,
\begin{equation}
\label{b4}
\pi'h=\frac{1}{2\pi}\int_{\Gamma^+}h(\xi)d\xi.
\end{equation}
So $\pi'(H^-)=0$. For $h\in H\bigcap L^1({\bf R})$, $\pi'h=\frac{1}{2\pi}\int_{{\bf R}}h(v)dv$ and for $h\in H^+\bigcap L^1({\bf R})$, $\pi'h=0$.\\
\indent An operator of order $m\in {\bf Z}$ and type $d$ is a matrix\\
$$\widetilde{A}=\left(\begin{array}{lcr}
  \pi^+P+G  & K  \\
   T  &  \widetilde{S}
\end{array}\right):
\begin{array}{cc}
\   C^{\infty}(M,E_1)\\
 \   \bigoplus\\
 \   C^{\infty}(\partial{M},F_1)
\end{array}
\longrightarrow
\begin{array}{cc}
\   C^{\infty}(M,E_2)\\
\   \bigoplus\\
 \   C^{\infty}(\partial{M},F_2)
\end{array},
$$
where $M$ is a manifold with boundary $\partial M$ and
$E_1,E_2$~ (respectively $F_1,F_2$) are vector bundles over $M~$ (respectively $\partial M
$).~Here,~$P:C^{\infty}_0(\Omega,\overline {E_1})\rightarrow
C^{\infty}(\Omega,\overline {E_2})$ is a classical
pseudodifferential operator of order $m$ on $\Omega$, where
$\Omega$ is a collar neighborhood of $M$ and
$\overline{E_i}|M=E_i~(i=1,2)$. $P$ has an extension:
$~{\cal{E'}}(\Omega,\overline {E_1})\rightarrow
{\cal{D'}}(\Omega,\overline {E_2})$, where
${\cal{E'}}(\Omega,\overline {E_1})~({\cal{D'}}(\Omega,\overline
{E_2}))$ is the dual space of $C^{\infty}(\Omega,\overline
{E_1})~(C^{\infty}_0(\Omega,\overline {E_2}))$. Let
$e^+:C^{\infty}(M,{E_1})\rightarrow{\cal{E'}}(\Omega,\overline
{E_1})$ denote extension by zero from $M$ to $\Omega$ and
$r^+:{\cal{D'}}(\Omega,\overline{E_2})\rightarrow
{\cal{D'}}(\Omega, {E_2})$ denote the restriction from $\Omega$ to
$X$, then define
$$\pi^+P=r^+Pe^+:C^{\infty}(M,{E_1})\rightarrow {\cal{D'}}(\Omega,
{E_2}).$$ In addition, $P$ is supposed to have the
transmission property; this means that, for all $j,k,\alpha$, the
homogeneous component $p_j$ of order $j$ in the asymptotic
expansion of the
symbol $p$ of $P$ in local coordinates near the boundary satisfies:\\
$$\partial^k_{x_n}\partial^\alpha_{\xi'}p_j(x',0,0,+1)=
(-1)^{j-|\alpha|}\partial^k_{x_n}\partial^\alpha_{\xi'}p_j(x',0,0,-1),$$
then $\pi^+P:C^{\infty}(M,{E_1})\rightarrow C^{\infty}(M,{E_2})$. Let $G$, $T$ be respectively the singular Green operator
and the {\rm trace} operator of order $m$ and type $d$. Let $K$ be a
potential operator and $S$ be a classical pseudodifferential
operator of order $m$ along the boundary. Denote by $B^{m,d}$ the collection of all operators of
order $m$
and type $d$,  and $\mathcal{B}$ is the union over all $m$ and $d$.\\
\indent Recall that $B^{m,d}$ is a Fr\'{e}chet space. The composition
of the above operator matrices yields a continuous map:
$B^{m,d}\times B^{m',d'}\rightarrow B^{m+m',{\rm max}\{
m'+d,d'\}}.$ Write $$\widetilde{A}=\left(\begin{array}{lcr}
 \pi^+P+G  & K \\
 T  &  \widetilde{S}
\end{array}\right)
\in B^{m,d},
 \widetilde{A}'=\left(\begin{array}{lcr}
\pi^+P'+G'  & K'  \\
 T'  &  \widetilde{S}'
\end{array} \right)
\in B^{m',d'}.$$\\
The composition $\widetilde{A}\widetilde{A}'$ is obtained by multiplication of the matrices (For more details see \cite{Ka}).
For
example $\pi^+P\circ G'$ and $G\circ G'$ are singular Green
operators of type $d'$ and
$$\pi^+P\circ\pi^+P'=\pi^+(PP')+L(P,P').$$
Here $PP'$ is the usual
composition of pseudodifferential operators and $L(P,P')$ called
leftover term is a singular Green operator of type $m'+d$. For our case, $P,~P'$ are classical pseudo differential operators, in other words $\pi^+P\in \mathcal{B}^{\infty}$ and $\pi^+P'\in \mathcal{B}^{\infty}$ .\\
\indent Let $M$ be a $n$-dimensional compact oriented manifold with boundary $\partial M$.
Denote by $\mathcal{B}$ the Boutet de Monvel's algebra. We recall that the main theorem in \cite{FGLS,Wa3}.
\begin{thm}\label{thm2.1}{\rm\cite{FGLS}} {\bf(Fedosov-Golse-Leichtnam-Schrohe)}
 Let $M$ and $\partial M$ be connected, ${\rm dim}M=n\geq3$, and let $\widetilde{S}$ (respectively $\widetilde{S}'$) be the unit sphere about $\xi$ (respectively $\xi'$) and $\sigma(\xi)$ (respectively $\sigma(\xi')$) be the corresponding canonical
$n-1$ (respectively $(n-2)$) volume form.
 Set $\widetilde{A}=\left(\begin{array}{lcr}\pi^+P+G &   K \\
T &  \widetilde{S}    \end{array}\right)$ $\in \mathcal{B}$ , and denote by $p$, $b$ and $s$ the local symbols of $P,G$ and $\widetilde{S}$ respectively.
 Define:
 \begin{align}\label{b5}
{\rm{\widetilde{Wres}}}(\widetilde{A})&=\int_X\int_{\bf \widetilde{ S}}{\rm{{\rm trace}}}_E\left[p_{-n}(x,\xi)\right]\sigma(\xi)dx \nonumber\\
&+2\pi\int_ {\partial X}\int_{\bf \widetilde{S}'}\left\{{\rm {\rm trace}}_E\left[({\rm{trace}}b_{-n})(x',\xi')\right]+{\rm{{\rm trace}}}
_F\left[s_{1-n}(x',\xi')\right]\right\}\sigma(\xi')dx',
\end{align}
where ${\rm{\widetilde{Wres}}}$ denotes the noncommutative residue of an operator in the Boutet de Monvel's algebra.\\
Then~~ a) ${\rm \widetilde{Wres}}([\widetilde{A},B])=0 $, for any
$\widetilde{A},B\in\mathcal{B}$;~~ b) It is the unique continuous {\rm trace} on
$\mathcal{B}/\mathcal{B}^{-\infty}$.
\end{thm}

\section{Spectral Einstein functionals for perturbations of Dirac operators}
\label{section:2.1}
Firstly, we recall the definition of the Dirac operator. Let $M$ be an $n$ dimensional oriented compact spin Riemannian manifold with a Riemannian metric $g^{M}$ and let $\nabla^L$ be the Levi-Civita connection about $g^{M}$.

In 
the fixed orthonormal frame $\{e_1,\cdots,e_n\}$, the connection matrix $(\omega_{s,t})$ is defined by
\begin{equation}
\label{a2}
\nabla^L(e_1,\cdots,e_n)= (e_1,\cdots,e_n)(\omega_{s,t}).
\end{equation}
\indent Let $c(X)$ be a Clifford action on $M$, where $X$ is a smooth vector field on $M$, which satisfies
\begin{align}
\label{a4}
&c(e_i)c(e_j)+c(e_j)c(e_i)=-2g^{M}(e_i,e_j).
\end{align}
In \cite{Y}, the Dirac operator is given
\begin{align}
\label{a5}
D=\sum^m_{i=1}c(e_i)\bigg[e_i-\frac{1}{4}\sum_{s,t}\omega_{s,t}
(e_i)
c(e_s)c(e_t)\bigg].
\end{align}
Set $$X=\sum\limits_{\alpha=1}^na_{\alpha}e_\alpha=X^T+X_n\partial_{x_n}=\sum\limits_{j=1}^nX_j\partial_j,$$ $$L(X):=\frac{1}{4}\sum\limits_{ij}\langle\nabla_X^L{e_i},e_j\rangle c(e_i)c(e_j)$$ and $$G(X,\Psi):=-\frac{1}{2}[c(X)c(\Psi)+c(\Psi) c(X)].$$
Let $\nabla_X^{S(TM)}$ denotes a spin connection, defined as
and
$$\nabla_X^{S(TM)}:=X+L(X).$$
And let $\nabla^{\Psi}_{X}=\nabla^{S(TM)}_{X}+G(X,\Psi)$ and $\nabla^{\Psi}_{Y}=\nabla^{S(TM)}_{Y}+G(Y,\Psi)$,
$c(\Psi)=c(X_{1})\cdots c(X_{k})$, $X_{i}~(1\leq i\leq k)$ are tangent vector fileds, $g^{ij}=g(dx_{i},dx_{j})$, $\xi=\sum\limits_{k}\xi_{j}dx_{j}$ and $\nabla^L_{\partial_{i}}\partial_{j}=\sum\limits_{k}\Gamma_{ij}^{k}\partial_{k}$,  we denote that
\begin{align}\label{a6}
&\sigma_{i}=-\frac{1}{4}\sum_{s,t}\omega_{s,t}
(e_i)c(e_i)c(e_s)c(e_t)
;~~~\xi^{j}=g^{ij}\xi_{i};~~~~\Gamma^{k}=g^{ij}\Gamma_{ij}^{k};~~~~\sigma^{j}=g^{ij}\sigma_{i}.
\end{align}
Then, perturbations of Dirac Operators is defined as
\begin{align}\label{a7}
D_{\Psi}=D+c(\Psi)=\sum^n_{i=1}c(e_i)\bigg[e_i-\frac{1}{4}\sum_{s,t}\omega_{s,t}
(e_i)
c(e_s)c(e_t)\bigg]+c(\Psi),
\end{align}
and
\begin{align}\label{a8}
\sigma_{1}(D_{\Psi})=&ic(\xi);\\
\sigma_{0}(D_{\Psi})=&-\frac{1}{4}\sum_{i,s,t}\omega_{s,t}(e_{i})c(e_{s})c(e_{t})
+c(\Psi),
\end{align}
where $\sigma_{l}(A)$ denotes the $l$-th order symbol of an operator $A$.

In addition, we obtain
\begin{align}\label{1000}
\nabla^{\Psi}_{X}\nabla^{\Psi}_{Y}&=\Big(X+L(X)
+G(X,\Psi)\Big)\Big(Y+L(Y)+G(Y,\Psi)\Big)\nonumber\\
 &=XY+X[L(Y)]+L(Y)X+L(X)Y+L(X)L(Y)+G(X,\Psi)Y
 +G(X,\Psi)L(Y)\nonumber\\
 &~~~~+X[G(Y,\Psi)]+G(Y,\Psi)X+L(X)G(Y,\Psi)+G(X,\Psi)G(Y,\Psi),
\end{align}
where
$X=\sum\limits_{j=1}^nX_j\partial_{x_j}, Y=\sum\limits_{l=1}^nY_l\partial_{x_l}$.

Let $p_{1},p_{2}$ be nonnegative integers and $p_{1}+p_{2}\leq n$,
 an application of (3.5) and (3.6) in \cite{Wa1} shows that
\begin{defn} \label{defn2.3}
Spectral Einstein functional of spin manifolds with boundary for perturbations of the Dirac operator are defined by
\begin{equation}\label{}
   \widetilde{Wres}[\pi^{+}\nabla^{\Psi}_{X}\nabla^{\Psi}_{Y}D_{\Psi}^{- p_{1}}
    \circ\pi^{+}D_{\Psi}^{-p_{2}}],
\end{equation}
where $\pi^{+}\nabla^{\Psi}_{X}\nabla^{\Psi}_{Y}D_{\Psi}^{- p_{1}}$, $\pi^{+}D_{\Psi}^{-p_{2}}$ are
 elements in Boutet de Monvel's algebra\cite{Wa3}.
\end{defn}
The following lemma derived from Dabrowski et al.'s Einstein functional plays a crucial role in our demonstration of the generalized noncommutative residue for the spectral action for perturbations of Dirac operators on compact manifolds with boundary.

 \begin{lem}\cite{DL}
The Einstein functional equal to
 \begin{equation*}
Wres\big(\widetilde{\nabla}_{X}\widetilde{\nabla}_{Y}\Delta^{-\frac{n}{2}}_{\Psi}\big)
=\frac{\upsilon_{n-1}}{6}2^{\frac{n}{2}}\int_{M}EG(X,Y)vol_{g}
 +\frac{\upsilon_{n-1}}{2}\int_{M}F(X,Y)vol_{g}+\frac{1}{2}\int_{M}({\rm trace}E)g(X,Y)vol_{g},
\end{equation*}
where $EG(X,Y)$ denotes the Einstein tensor evaluated on the two vector fields, $F(X,Y)=Tr(X_{a}Y_{b}F_{ab})$ and
$F_{ab}$ is the curvature tensor of the connection $T$, $\mathrm{Tr}E$ denotes the trace of $E$ and $\upsilon_{n-1}=\frac{2\pi^{\frac{n}{2}}}{\Gamma(\frac{n}{2})}$.
\end{lem}

Consider $X$ and $Y$ as a pair of vector fields on a compact Riemannian manifold $M$ with dimension $n$. By employing the Laplace operator $\Delta_{\Psi}=D_{\Psi}^{2}=\Delta+E$, which acts on sections of a vector bundle $S(TM)$ with a rank of $2^{\frac{n}{2}}$.
The purpose of this section is to demonstrate the following.
\begin{thm}\label{theorem10000}
For the Laplace (type) operator with torsion $\Delta_{\Psi}$, the Einstein functional equal to
\begin{align*}
Wres\big(\nabla^{\Psi}_{X}\nabla^{\Psi}_{Y}\Delta^{-\frac{n}{2}}_{\Psi}\big)
=&\frac{\upsilon_{n-1}}{6}2^{\frac{n}{2}}\int_{M}EG(X,Y)vol_{g}
+\frac{\upsilon_{n-1}}{4}\int_{M}{\rm trace}\Big[c(X)\nabla^{S(TM)}_Y\big(c(\Psi)\big)\nonumber\\
&+\nabla^{S(TM)}_Y\big(c(\Psi)\big) c(X)-c(Y)\nabla^{S(TM)}_X(c(\Psi))-\nabla^{S(TM)}_X(c(\Psi)) c(Y)\Big]dvol_{M}\nonumber\\
&+\frac{1}{2}\int_{M}{\rm trace}[\frac{1}{4}s+\sum\limits^{n}_{j=1}\frac{1}{2}c(\Psi)c(e_j)
c(\Psi)c(e_j)+(1-\frac{n}{2})\big(c(\Psi)\big)^2]g(X,Y)dvol_{M},
\end{align*}
where 
$s$ is the scalar curvature.
\begin{proof}
Let $X=\sum\limits_{a=1}^{n}X^{a}e_{a}$, $Y=\sum\limits_{b=1}^{n}Y^{b}e_{b}$.
Considering that
 \begin{equation}
F(X,Y)={\rm trace}(X_{a}Y_{b}F_{ab})=\sum_{a,b=1}^{n}X^{a}Y^{b}{\rm trace}^{S(TM)}(F_{e_{a},e_{b}}).
\end{equation}
Set $\overline{A}(X)=L(X)+G(X,\Psi)$,
then we acquire
 \begin{align}
F_{e_{a},e_{b}}=&(e_{a}+\overline{A}(e_{a}))(e_{b}+\overline{A}(e_{b}))-(e_{b}+\overline{A}(e_{b}))(e_{a}+\overline{A}(e_{a}))
-([e_{a},e_{a}]+\overline{A}([e_{a},e_{b}]))
\nonumber\\
=&e_{a}\circ \overline{A}(e_{b})+\overline{A}(e_{a})\circ e_{b}+\overline{A}(e_{a})A(e_{b})-e_{b}\circ \overline{A}(e_{a})
-\overline{A}(e_{b})\circ e_{a}\nonumber\\
& -\overline{A}(e_{b})\overline{A}(e_{a})-\overline{A}([e_{a},e_{b}])\nonumber\\
=&\overline{A}(e_{b})e_{a}+e_{a}(\overline{A}(e_{b}))+\overline{A}(e_{a})\circ e_{b}+\overline{A}(e_{a})\overline{A}(e_{b})
-\overline{A}(e_{a})\circ e_{b}-e_{b}(\overline{A}(e_{a}))\nonumber\\
& -\overline{A}(e_{b})e_{a}-\overline{A}(e_{b})\overline{A}(e_{a})-\overline{A}([e_{a},e_{b}])\nonumber\\
=&e_{a}(\overline{A}(e_{b}))-e_{b}(\overline{A}(e_{a}))+\overline{A}(e_{a})\overline{A}(e_{b})
-\overline{A}(e_{b})\overline{A}(e_{a})-\overline{A}([e_{a},e_{b}]).
\end{align}
We note that ${\rm trace}[\overline{A}(e_{a})\overline{A}(e_{b})-\overline{A}(e_{b})\overline{A}(e_{a})]=0 $ and $e_{a}(c(e_s))=0$. If $s=t$, we have $\omega_{s,t}(e_{b})=0$
; if $s\neq t$, we have ${\rm trace}[c(e_s)c(e_t)]=0$. Then we obtain ${\rm trace}[e_{a}\big(\omega_{s,t}(e_{b}) \big)c(e_s)c(e_t)]=0$,
so there is the following formula
 \begin{align}
{\rm trace}\big(e_{a}(\overline{A}(e_{b}))\big)
=&{\rm trace}\Big[e_{a}\Big(-\frac{1}{4}\sum_{s,t} \omega_{s,t}(e_{b}) c(e_s)c(e_t)
-\frac{1}{2}[c(e_{b})c(\Psi)+c(\Psi) c(e_{b})]\Big)\Big]\nonumber\\
=&-\frac{1}{2}{\rm trace}\Big[c(e_{b})e_{a}\Big(c(\Psi) \Big)+e_{a}\Big(c(\Psi) \Big) c(e_{b})\Big]
\end{align}
and
\begin{align}
&{\rm trace}\big(\overline{A}([e_{a},e_{b}])\big)\nonumber\\
=&{\rm trace} \Big(-\frac{1}{4}\sum_{s,t} \omega_{s,t}([e_{a},e_{b}]) c(e_s)c(e_t)-\frac{1}{2}[c([e_{a},e_{b}])c(\Psi)+c(\Psi) c([e_{a},e_{b}])\Big)(x_{0})\nonumber\\
=&-\frac{1}{2}{\rm trace} \Big(c([e_{a},e_{b}])c(\Psi)+c(\Psi) c([e_{a},e_{b}])\Big).
\end{align}
Fix a point $x_0$, $\omega_{st}(e_a)(x_0)=0$, then $e_a=\nabla^{S(TM)}_{e_a}(x_0)$ and let $e_a=\sum\limits_{j}H_{a_j}\partial_{j}$, then we have
\begin{align}
[e_a,e_b](x_0)&=[H_{a_j}\partial_{j},H_{b_l}\partial_{l}](x_0)\nonumber\\
&=H_{a_j}\partial_{j}(H_{b_l})\partial_{l}(x_0)-H_{b_l}\partial_{l}(H_{a_j})\partial_{j}(x_0)\nonumber\\
&=0
\end{align}
we therefore draw the following formula
\begin{align}
&\sum\limits^{n}_{a,b=1}X^{a}Y^{b}\bigg\{-\frac{1}{2}{\rm trace}\Big[c(e_{b})e_{a}\Big(c(\Psi) \Big)+e_{a}\Big(c(\Psi) \Big) c(e_{b})\Big]+\frac{1}{2}{\rm trace}\Big[c(e_{a})e_{b}\Big(c(\Psi) \Big)\nonumber\\
&+e_{b}\Big(c(\Psi) \Big) c(e_{a})\Big]+\frac{1}{2}{\rm trace} \Big(c([e_{a},e_{b}])c(\Psi)+c(\Psi) c([e_{a},e_{b}])\Big)  \bigg\}\nonumber\\
=&-\frac{1}{2}{\rm trace}\Big(c(Y)\nabla^{S(TM)}_X(c(\Psi))+\nabla^{S(TM)}_X(c(\Psi)) c(Y)\Big)\nonumber\\
&+\frac{1}{2}{\rm trace}\Big(c(X)\nabla^{S(TM)}_Y\big(c(\Psi)\big)+\nabla^{S(TM)}_Y\big(c(\Psi)\big) c(X)\Big)\nonumber\\
=&F(X,Y).
\end{align}
Let $\Delta_{\Psi}=\Delta+E$, as the formula (2,28) in \cite{AAA}, we have
\begin{align}
{\rm trace}(E)
=&{\rm trace}[\frac{1}{4}s+\frac{1}{2}\Psi c(e_j)
\Psi c(e_j)+(1-\frac{n}{2})\big(c(\Psi)\big)^2]\nonumber\\
=&{\rm trace}[\frac{1}{4}s+\sum\limits^{n}_{j=1}\frac{1}{2}c(\Psi)c(e_j)
c(\Psi)c(e_j)+(1-\frac{n}{2})\big(c(\Psi)\big)^2].
\end{align}
Taking all the above conclusions into consideration, we ultimately complete the proof of the Theorem.
\end{proof}

\end{thm}
\section{ The noncommutative residue $\widetilde{{\rm Wres}}\bigg[\pi^{+}\nabla^{\Psi}_{X}\nabla^{\Psi}_{Y}D_{\Psi}^{-2}
    \circ\pi^{+}D_{\Psi}^{-2}\bigg]$ on manifolds with boundary}
\label{section:3}
In this section, we calculate the spectral Einstein functional for 4-dimension compact manifolds with boundary and derive a Dabrowski-Sitarz-Zalecki type formula in this case.

For perturbations of Dirac Operators
 $\nabla^{\Psi}_{X}\nabla^{\Psi}_{Y}D_{\Psi}^{-2}$ and $D_{\Psi}^{-2}$,
 let $\sigma_{l}(A)$ denote the $l$-order symbol of an operator A. Then, based on $\sigma(\partial_{x_j})=i\xi_j$ and (\ref{1000}), we can establish the following lemmas.
\begin{lem}\label{lem2} The following identities hold:
\begin{align}\label{b22}
\sigma_{0}(\nabla_X^{\Psi}\nabla_Y^{\Psi})=&X[L(Y)]+L(X)L(Y)+X[G(Y,\Psi)]+G(X,\Psi)L(Y)+L(X)G(Y,\Psi)\nonumber\\
&+G(X,\Psi)G(Y,\Psi);\nonumber\\
\sigma_{1}(\nabla_X^{\Psi}\nabla_Y^{\Psi})=&i\sum_{j,l=1}^nX_j\frac{\partial_{Y_l}}{\partial_{x_j}}i\xi_l
+i\sum_jA(Y)X_j\xi_j+i\sum_lA(Y)Y_l\xi_l+\sum_j  G(X,\Psi)Y_ji\xi_j\nonumber\\
&+\sum_j G(X,\Psi)X_ji \xi_j;\nonumber\\
\sigma_{2}(\nabla_X^{\Psi}\nabla_Y^{\Psi})=&-\sum\limits_{j,l=1}^nX_jY_l\xi_j\xi_l.
\end{align}
\end{lem}

Next, we present the following lemmas.
\begin{lem}\label{lem3} The following identities hold:
\begin{align}\label{4.8}
\sigma_{-1}(D_{\Psi}^{-1})&=\frac{ic(\xi)}{|\xi|^{2}}; \\
\sigma_{-2}(D_{\Psi}^{-1})&=\frac{c(\xi)\sigma_{0}(D_{\Psi})c(\xi)}{|\xi|^{4}}+\frac{c(\xi)}{|\xi|^{6}}\sum_{j}c(\texttt{d}x_{j})
\Big[\partial_{x_{j}}(c(\xi))|\xi|^{2}-c(\xi)\partial_{x_{j}}(|\xi|^{2})\Big]; \\
\sigma_{-2}(D_{\Psi}^{-2})&=\sigma_{-2}(D^{-2})=|\xi|^{-2};\\
\sigma_{-3}(D_{\Psi}^{-2})&=-i|\xi|^{-4}\xi_k(\Gamma^k-2\delta^k)
-i|\xi|^{-6}2\xi^j\xi_\alpha\xi_\beta\partial_jg^{\alpha\beta}-\bigg(c(\Psi)ic(\xi)+ic(\xi)c(\Psi)  \bigg)|\xi|^{-4}.
\end{align}
\end{lem}
According to Lemma \ref{lem2} and Lemma \ref{lem3}, we obtain
\begin{lem} \label{lem100}
The following identities hold:
\begin{align}\label{4.10}
\sigma_{0}(\nabla^{\Psi}_{X}\nabla^{\Psi}_{Y}D_{\Psi}^{-2})=&
-\sum_{j,l=1}^nX_jY_l\xi_j\xi_l|\xi|^{-2};\\
\sigma_{-1}(\nabla^{\Psi}_{X}\nabla^{\Psi}_{Y}D_{\Psi}^{-2})=&
\sigma_{2}(\nabla^{\Psi}_{X}\nabla^{\Psi}_{Y})\sigma_{-3}(D_{\Psi}^{-2})
+\sigma_{1}(\nabla^{\Psi}_{X}\nabla^{\Psi}_{Y})\sigma_{-2}(D_{\Psi}^{-2})\nonumber\\
&+\sum_{j=1}^{n}\partial_{\xi_{j}}\big[\sigma_{2}(\nabla^{\Psi}_{X}\nabla^{\Psi}_{Y})\big]
D_{x_{j}}\big[\sigma_{-2}(D_{\Psi}^{-2})\big].
\end{align}
\end{lem}
Since $\Theta$ is a global form on $\partial M$, so for any fixed point $x_{0}\in\partial M$, we can choose the normal coordinates
$U$ of $x_{0}$ in $\partial M$(not in $M$) and compute $\Theta(x_{0})$ in the coordinates $\widetilde{U}=U\times [0,1)$ and the metric
$\frac{1}{h(x_{n})}g^{\partial M}+dx _{n}^{2}$. The dual metric of $g^{M}$ on $\widetilde{U}$ is
$h(x_{n})g^{\partial M}+dx _{n}^{2}.$ Write
$g_{ij}^{M}=g^{M}(\frac{\partial}{\partial x_{i}},\frac{\partial}{\partial x_{j}})$;
$g^{ij}_{M}=g^{M}(d x_{i},dx_{j})$, then

\begin{equation*}
[g_{i,j}^{M}]=
\begin{bmatrix}\frac{1}{h( x_{n})}[g_{i,j}^{\partial M}]&0\\0&1\end{bmatrix};\quad
[g^{i,j}_{M}]=\begin{bmatrix} h( x_{n})[g^{i,j}_{\partial M}]&0\\0&1\end{bmatrix},
\end{equation*}
and
\begin{equation*}
\partial_{x_{s}} g_{ij}^{\partial M}(x_{0})=0,\quad 1\leq i,j\leq n-1;\quad g_{i,j}^{M}(x_{0})=\delta_{ij}.
\end{equation*}

Let $\{e_{1},\cdots, e_{n-1}\}$ be an orthonormal frame field in $U$ about $g^{\partial M}$ which is parallel along geodesics and
$e_{i}=\frac{\partial}{\partial x_{i}}(x_{0})$, then $\{\widetilde{e_{1}}=\sqrt{h(x_{n})}e_{1}, \cdots,
\widetilde{e_{n-1}}=\sqrt{h(x_{n})}e_{n-1},
\widetilde{e_{n}}=dx_{n}\}$ is the orthonormal frame field in $\widetilde{U}$ about $g^{M}.$
Locally $S(TM)|\widetilde{U}\cong \widetilde{U}\times\wedge^{*}_{C}(\frac{n}{2}).$ Let $\{f_{1},\cdots,f_{n}\}$ be the orthonormal basis of
$\wedge^{*}_{C}(\frac{n}{2})$. Take a spin frame field $\sigma: \widetilde{U}\rightarrow Spin(M)$ such that
$\pi\sigma=\{\widetilde{e_{1}},\cdots, \widetilde{e_{n}}\}$ where $\pi: Spin(M)\rightarrow O(M)$ is a double covering, then
$\{[\sigma, f_{i}], 1\leq i\leq n\}$ is an orthonormal frame of $S(TM)|_{\widetilde{U}}$. In the following, since the global form $\Theta$
is independent of the choice of the local frame, so we can compute ${\rm {\rm trace}}_{S(TM)}$ in the frame $\{[\sigma, f_{i}], 1\leq i\leq n\}$.
Let $\{\hat{e}_{1},\cdots,\hat{e}_{n}\}$ be the canonical basis of $\mathbb{R}^{n}$ and
$c(\hat{e}_{i})\in cl_{C}(n)\cong Hom(\wedge^{*}_{C}(\frac{n}{2}),\wedge^{*}_{C}(\frac{n}{2}))$ be the Clifford action. Then
\begin{equation*}
c(\widetilde{e_{i}})=[(\sigma,c(\hat{e}_{i}))]; \quad c(\widetilde{e_{i}})[(\sigma, f_{i})]=[\sigma,(c(\hat{e}_{i}))f_{i}]; \quad
\frac{\partial}{\partial x_{i}}=[(\sigma,\frac{\partial}{\partial x_{i}})],
\end{equation*}
then we have $\frac{\partial}{\partial x_{i}}c(\widetilde{e_{i}})=0$ in the above frame. In accordance with Lemma 2.2 in \cite{Wa3}, we obtain

\begin{lem}\label{le:32}
With the metric $g^{M}$ on $M$ near the boundary
\begin{eqnarray}\label{fum100}
\partial_{x_j}(|\xi|_{g^M}^2)(x_0)&=&\left\{
       \begin{array}{c}
        0,  ~~~~~~~~~~ ~~~~~~~~~~ ~~~~~~~~~~~~~{\rm if }~j<n; \\[2pt]
       h'(0)|\xi'|^{2}_{g^{\partial M}},~~~~~~~~~~~~~~~~~~~~~{\rm if }~j=n.
       \end{array}
    \right. \\
\partial_{x_j}[c(\xi)](x_0)&=&\left\{
       \begin{array}{c}
      0,  ~~~~~~~~~~ ~~~~~~~~~~ ~~~~~~~~~~~~~{\rm if }~j<n;\\[2pt]
\partial x_{n}(c(\xi'))(x_{0}), ~~~~~~~~~~~~~~~~~{\rm if }~j=n,
       \end{array}
    \right.
\end{eqnarray}
where $\xi=\xi'+\xi_{n}dx_{n}$.
\end{lem}
 An application of (2.1.4) in \cite{Wa1} shows that
\begin{align}\label{fum100}
&\widetilde{{\rm Wres}}\bigg[\pi^{+}\nabla^{\Psi}_{X}\nabla^{\Psi}_{Y}D_{\Psi}^{-p_{1}}
    \circ\pi^{+}D_{\Psi}^{-p_{2}}\bigg]\nonumber\\
&=\int_{M}\int_{|\xi|=1}{\rm trace}_{S(TM)}
  [\sigma_{-n}(\nabla^{\Psi}_{X}\nabla^{\Psi}_{Y}D_{\Psi}^{-p_{1}}
  \circ D_{\Psi}^{-p_{2}}]\sigma(\xi)\texttt{d}x+\int_{\partial M}\Phi,
\end{align}
where
 \begin{align}\label{fum101}
\Phi=&\int_{|\xi'|=1}\int_{-\infty}^{+\infty}\sum_{j,k=0}^{\infty}\sum \frac{(-i)^{|\alpha|+j+k+\ell}}{\alpha!(j+k+1)!}
{\rm trace}_{S(TM)}[\partial_{x_{n}}^{j}\partial_{\xi'}^{\alpha}\partial_{\xi_{n}}^{k}\sigma_{r}^{+}
((\nabla^{\Psi}_{X}\nabla^{\Psi}_{Y}D_{\Psi}^{-p_{1}})(x',0,\xi',\xi_{n})\nonumber\\
&\times\partial_{x_{n}}^{\alpha}\partial_{\xi_{n}}^{j+1}\partial_{x_{n}}^{k}\sigma_{l}(D_{\Psi}^{-p_{2}})(x',0,\xi',\xi_{n})]
\texttt{d}\xi_{n}\sigma(\xi')\texttt{d}x' ,
\end{align}
and the sum is taken over $r-k+|\alpha|+\ell-j-1=-n,r\leq-p_{1},\ell\leq-p_{2}$.

In the following, we will compute the residue $\widetilde{{\rm Wres}}\bigg[\pi^{+}\nabla^{\Psi}_{X}\nabla^{\Psi}_{Y}D_{\Psi}^{-2}
    \circ\pi^{+}D_{\Psi}^{-2}\bigg]$ on 4-dimensional oriented
compact spin manifolds with boundary and get a  Dabrowski-Sitarz-Zalecki
type theorem in this case. By  (\ref{fum100}) and (\ref{fum101}), we have
\begin{align}\label{fum102}
&\widetilde{{\rm Wres}}\bigg[\pi^{+}\nabla^{\Psi}_{X}\nabla^{\Psi}_{Y}D_{\Psi}^{-2}
    \circ\pi^{+}D_{\Psi}^{-2}\bigg]\nonumber\\
&=\int_{M}\int_{|\xi|=1}{\rm {\rm trace}}_{S(TM)}
  [\sigma_{-4}(\nabla^{\Psi}_{X}\nabla^{\Psi}_{Y}D_{\Psi}^{-2}
  \circ D_{\Psi}^{-2}]\sigma(\xi)\texttt{d}x+\int_{\partial M}\widetilde{\Phi},
\end{align}
where
 \begin{align}\label{fum103}
\widetilde{\Phi}=&\int_{|\xi'|=1}\int_{-\infty}^{+\infty}\sum_{j,k=0}^{\infty}\sum \frac{(-i)^{|\alpha|+j+k+\ell}}{\alpha!(j+k+1)!}
{\rm {\rm trace}}_{S(TM)}[\partial_{x_{n}}^{j}\partial_{\xi'}^{\alpha}\partial_{\xi_{n}}^{k}\sigma_{r}^{+}
(\nabla^{\Psi}_{X}\nabla^{\Psi}_{Y}D_{\Psi}^{-2})(x',0,\xi',\xi_{n})\nonumber\\
&\times\partial_{x_{n}}^{\alpha}\partial_{\xi_{n}}^{j+1}\partial_{x_{n}}^{k}\sigma_{l}(D_{\Psi}^{-2})(x',0,\xi',\xi_{n})]
\texttt{d}\xi_{n}\sigma(\xi')\texttt{d}x' ,
\end{align}
and the sum is taken over $r-k+|\alpha|+\ell-j-1=-4,r\leq0,\ell\leq-2$.

Locally we can use Theorem \ref{theorem10000} to compute the interior term of (\ref{fum102}), then
 \begin{eqnarray}\label{fum3.10}
&&\int_{M}\int_{|\xi|=1}{\rm {\rm trace}}_{S(TM)}
  [\sigma_{-4}(\nabla^{\Psi}_{X}\nabla^{\Psi}_{Y}D_{\Psi}^{-2}
  \circ D_{\Psi}^{-2}]\sigma(\xi)\texttt{d}x\nonumber\\
&=&\frac{4\pi^2}{3}\int_{M}EG(X,Y)vol_{g}
+\frac{\pi^2}{2}\int_{M}{\rm trace}\Big[c(X)\nabla^{S(TM)}_Y\big(c(\Psi)\big)\nonumber\\
&&+\nabla^{S(TM)}_Y\big(c(\Psi)\big) c(X)-c(Y)\nabla^{S(TM)}_X(c(\Psi))-\nabla^{S(TM)}_X(c(\Psi)) c(Y)\Big]dvol_{M}\nonumber\\
&&+\frac{1}{2}\int_{M}{\rm trace}[\frac{1}{4}s+\sum\limits^{n}_{j=1}\frac{1}{2}c(\Psi)c(e_j)
c(\Psi)c(e_j)-\big(c(\Psi)\big)^2]g(X,Y)dvol_{M},
\end{eqnarray}
so we only need to compute $\int_{\partial M}\widetilde{\Phi}$.\\

 When $n=4$, then ${\rm {\rm trace}}_{S(TM)}[{\rm id}]=4$, the sum is taken over $r-k+|\alpha|+\ell-j=-3,r\leq0,\ell\leq-2$, then we have the $\int_{\partial{M}}\widetilde{\Phi}$
is the sum of the following five cases:
~\\

\noindent  {\bf case (a)~(I)}~$r=0,~l=-2,~k=j=0,~|\alpha|=1$\\

\noindent
By (\ref{fum103}), we get
\begin{equation}\label{b24}
\widetilde{\Phi}_1=-\int_{|\xi'|=1}\int^{+\infty}_{-\infty}\sum_{|\alpha|=1}
 {\rm trace}[\partial^\alpha_{\xi'}\pi^+_{\xi_n}\sigma_{0}(\nabla_X^{\Psi}\nabla_Y^{\Psi}D^{-2}_{\Psi})\times
 \partial^\alpha_{x'}\partial_{\xi_n}\sigma_{-2}(D^{-2}_{\Psi})](x_0)d\xi_n\sigma(\xi')dx'.
\end{equation}
By Lemma 2.2 in \cite{Wa3}, for $i<n$, then
\begin{equation}\label{b25}
\partial_{x_i}\sigma_{-2}(D^{-2}_{\Psi})(x_0)=
\partial_{x_i}(|\xi|^{-2})(x_0)=
-\frac{\partial_{x_i}(|\xi|^{2})(x_0)}{|\xi|^4}=0,
\end{equation}
 so $\widetilde{\Phi}_1=0$ and {\bf case (a)~(I)} vanishes.\\

 \noindent  {\bf case (a)~(II)}~$r=0,~l=-2,~k=|\alpha|=0,~j=1$.\\

\noindent By (\ref{fum103}), we get
\begin{equation}\label{200}
\widetilde{\Phi}_2=-\frac{1}{2}\int_{|\xi'|=1}\int^{+\infty}_{-\infty} {\rm
{\rm trace}} [\partial_{x_n}\pi^+_{\xi_n}\sigma_{0}(\nabla_X^{\Psi}\nabla_Y^{\Psi}D^{-2}_{\Psi})\times
\partial_{\xi_n}^2\sigma_{-2}(D^{-2}_{\Psi})](x_0)d\xi_n\sigma(\xi')dx'.
\end{equation}
By Lemma \ref{lem3}, we have
\begin{eqnarray}\label{201}
\partial_{\xi_n}^2\sigma_{-2}(D^{-2}_{\Psi})(x_0)
=\partial_{\xi_n}^2(|\xi|^{-2})(x_0)=\frac{6\xi_n^2-2}{(1+\xi_n^2)^3}.
\end{eqnarray}
It follows that
\begin{eqnarray}\label{202}
\partial_{x_n}\sigma_{0}(\nabla^{\Psi}_{X}\nabla^{\Psi}_{Y}D_{\Psi}^{-2})(x_0)
=\partial_{x_n}(-\sum\limits_{j,l=1}^nX_jY_l\xi_j\xi_l|\xi|^{-2})=\frac{\sum\limits_{j,l=1}^nX_jY_l\xi_j\xi_lh'(0)|\xi'|^2}{(1+\xi_n^2)^2}.
\end{eqnarray}
By integrating formula, we obtain
\begin{align}\label{203}
&\pi^+_{\xi_n}\partial_{x_n}\sigma_{0}(\nabla^{\Psi}_{X}\nabla^{\Psi}_{Y}D_{\Psi}^{-2})(x_0)\nonumber\\
=&\partial_{x_n}\pi^+_{\xi_n}\sigma_{0}(\nabla^{\Psi}_{X}\nabla^{\Psi}_{Y}D_{\Psi}^{-2})\nonumber\\
=&-\frac{i\xi_n}{4(\xi_n-i)^2}\sum\limits_{j,l=1}^{n-1}X_jY_l\xi_j\xi_lh'(0)+\frac{2-i\xi_n}{4(\xi_n-i)^2}X_nY_nh'(0)-\frac{i}{4(\xi_n-i)^2}\sum\limits_{j=1}^{n-1}X_jY_n\xi_j\nonumber\\
&-\frac{i}{4(\xi_n-i)^2}\sum\limits_{l=1}^{n-1}X_nY_l\xi_l.
\end{align}
We note that $i<n,~\int_{|\xi'|=1}\xi_{i_{1}}\xi_{i_{2}}\cdots\xi_{i_{2d+1}}\sigma(\xi')=0$,
so we omit some items that have no contribution for computing {\bf case (a)~(II)}.
From (\ref{201}) and (\ref{203}), we obtain
\begin{align}\label{204}
&{\rm {\rm trace}} [\partial_{x_n}\pi^+_{\xi_n}\sigma_{0}(\nabla^{\Psi}_{X}\nabla^{\Psi}_{Y}D_{\Psi}^{-2})\times
\partial_{\xi_n}^2\sigma_{-2}(D_{\Psi}^{-2})](x_0)\nonumber\\
=&\frac{2+2\xi_ni-6\xi_n^3i-2i}{(\xi_n-i)^5(\xi_n+i)^3}\sum\limits_{j,l=1}^{n-1}X_jY_l\xi_j\xi_lh'(0)
+\frac{2+2\xi_ni-6\xi_n^3i-2i}{(\xi_n-i)^5(\xi_n+i)^3}X_nY_nh'(0)\nonumber\\
&+\frac{2(1-3\xi_n^2)i}{(\xi_n-i)^5(\xi_n+i)^3}\sum\limits_{j=1}^{n-1}X_jY_n\xi_j
+\frac{2(1-3\xi_n^2)i}{(\xi_n-i)^5(\xi_n+i)^3}\sum\limits_{l=1}^{n-1}X_nY_l\xi_l.
\end{align}
Therefore, we get
\begin{align}\label{205}
\widetilde{\Phi}_2
=&-\frac{1}{2}\int_{|\xi'|=1}\int^{+\infty}_{-\infty}\bigg\{\frac{2+2\xi_ni-6\xi_n^3i-2i}
{(\xi_n-i)^5(\xi_n+i)^3}\sum\limits_{j,l=1}^{n-1}X_jY_l\xi_j\xi_lh'(0)
+\frac{2+2\xi_ni-6\xi_n^3i-2i}{(\xi_n-i)^5(\xi_n+i)^3}X_nY_nh'(0)\bigg\}\nonumber\\
&\times d\xi_n\sigma(\xi')dx'\nonumber\\
=&-\sum\limits_{j,l=1}^{n-1}X_jY_lh'(0)\Omega_3\int_{\Gamma^{+}}\frac{1+\xi_ni-3\xi_n^3i-i}
{(\xi_n-i)^5(\xi_n+i)^3}\xi_j\xi_ld\xi_{n}dx'-X_nY_nh'(0)\Omega_3\int_{\Gamma^{+}}
\frac{1+\xi_ni-3\xi_n^3i-i}{(\xi_n-i)^5(\xi_n+i)^3}\nonumber\\
&\times d\xi_{n}dx'\nonumber\\
=&-\frac{2\pi i}{4!}h'(0)\bigg(\sum\limits_{j,l=1}^{n-1}X_jY_l\left[\frac{1+\xi_ni-3\xi_n^3i-i}
{(\xi_n+i)^3}\right]^{(4)}\bigg|_{\xi_n=i}+X_nY_n
\left[\frac{1+\xi_ni-3\xi_n^3i-i}{(\xi_n+i)^3}\right]^{(4)}\bigg|_{\xi_n=i}\bigg)\Omega_3dx'\nonumber\\
=&\frac{13}{8}h'(0)\pi\Omega_3dx'\left(\frac{\pi}{3}\sum\limits_{j=1}^{n-1}X_jY_j+\frac{1}{4}X_nY_n\right),
\end{align}
where ${\rm \Omega_{3}}$ is the canonical volume of $S^{2}.$\\

\noindent  {\bf case (a)~(III)}~$r=0,~l=-2,~j=|\alpha|=0,~k=1.$\\

\noindent By (\ref{fum103}), we get
\begin{align}\label{300}
\widetilde{\Phi}_3&=-\frac{1}{2}\int_{|\xi'|=1}\int^{+\infty}_{-\infty}
{\rm {\rm trace}} [\partial_{\xi_n}\pi^+_{\xi_n}\sigma_{0}(\nabla_X^{\Psi}\nabla_Y^{\Psi}D^{-2}_{\Psi})\times
\partial_{\xi_n}\partial_{x_n}\sigma_{-2}(D^{-2}_{\Psi})](x_0)d\xi_n\sigma(\xi')dx'\nonumber\\
&=\frac{1}{2}\int_{|\xi'|=1}\int^{+\infty}_{-\infty}
{\rm {\rm trace}} [\partial_{\xi_n}^2\pi^+_{\xi_n}\sigma_{0}(\nabla_X^{\Psi}\nabla_Y^{\Psi}D^{-2}_{\Psi})\times
\partial_{x_n}\sigma_{-2}(D^{-2}_{\Psi})](x_0)d\xi_n\sigma(\xi')dx'.
\end{align}
Based on Lemma \ref{lem3}, we can express the equation as follows:
\begin{eqnarray}\label{b301}
\partial_{x_n}\sigma_{-2}(D_{\Psi}^{-2})(x_0)|_{|\xi'|=1}
=-\frac{h'(0)}{(1+\xi_n^2)^2}.
\end{eqnarray}
A straightforward computation yields
\begin{align}\label{b302}
\pi^+_{\xi_n}\sigma_{0}(\nabla^{\Psi}_{X}\nabla^{\Psi}_{Y}D^{-2}_{\Psi})(x_0)|_{|\xi'|=1}
=&\frac{i}{2(\xi_n-i)}\sum_{j,l=1}^{n-1}X_jY_l\xi_j\xi_l-\frac{1}{2(\xi_n-i)}X_nY_n\nonumber\\
&-\frac{1}{2(\xi_n-i)}\sum_{j=1}^{n-1}X_jY_n\xi_j-\frac{1}{2(\xi_n-i)}\sum_{l=1}^{n-1}X_nY_l\xi_l.
\end{align}
Additionally, simple calculations provide
\begin{align}\label{b303}
\partial_{\xi_n}^2\pi^+_{\xi_n}\sigma_{0}(\nabla^{\Psi}_{X}\nabla^{\Psi}_{Y}D^{-2}_{\Psi})(x_0)|_{|\xi'|=1}
=\frac{i}{(\xi_n-i)^3}\sum_{j,l=1}^{n-1}X_jY_l\xi_j\xi_l-\frac{1}{(\xi_n-i)^3}X_nY_n.
\end{align}
From (\ref{b301}) and (\ref{b303}), we obtain
\begin{align}\label{b30430}
&{\rm trace}[\partial_{\xi_n}\pi^+_{\xi_n}\sigma_{0}(\nabla^{\Psi}_{X}\nabla^{\Psi}_{Y}D^{-2}_{\Psi})\times
\partial_{\xi_n}\partial_{x_n}\sigma_{-2}(D_{\Psi}^{-2})](x_0)\nonumber\\
=&-4\frac{h'(0)i}{(\xi_n-i)^5(\xi_n+i)^2}\sum_{j,l=1}^{n-1}X_jY_l\xi_j\xi_l+4\frac{h'(0)}{(\xi_n-i)^5(\xi_n+i)^2}X_nY_n.
\end{align}
Hence, we obtain
\begin{align}\label{b30431}
\widetilde{\Phi}_3&=\frac{1}{2}\int_{|\xi'|=1}\int^{+\infty}_{-\infty}
\bigg(-4\frac{h'(0)i}{(\xi_n-i)^5(\xi_n+i)^2}
\sum_{j,l=1}^{n-1}X_jY_l\xi_j\xi_l+4\frac{h'(0)}{(\xi_n-i)^5(\xi_n+i)^2}X_nY_n\bigg)d\xi_n\sigma(\xi')dx'\nonumber\\
&=-2\sum_{j,l=1}^{n-1}X_jY_lh'(0)\Omega_3\int_{\Gamma^{+}}\frac{i}{(\xi_n-i)^5(\xi_n+i)^2}\xi_j\xi_ld\xi_{n}dx'+2X_nY_nh'(0)\Omega_3\int_{\Gamma^{+}}\frac{1}{(\xi_n-i)^5(\xi_n+i)^2}d\xi_{n}dx'\nonumber\\
&=-2\sum_{j,l=1}^{n-1}X_jY_lh'(0)\Omega_3\frac{2\pi i}{4!}\left[\frac{i}{(\xi_n+i)^2}\right]^{(4)}
\bigg|_{\xi_n=i}dx'+2X_nY_nh'(0)\Omega_3\frac{2\pi i}{4!}\left[\frac{1}{(\xi_n+i)^2}\right]^{(4)}\bigg|_{\xi_n=i}dx'\nonumber\\
&=\frac{5}{4}h'(0)\pi\Omega_3dx'\left(\frac{\pi}{3}\sum_{j=1}^{n-1}X_jY_j+\frac{i}{4}X_nY_n\right).
\end{align}

\noindent  {\bf case (b)}~$r=0,~l=-3,~k=j=|\alpha|=0$.\\

\noindent According to equation (\ref{fum103}), we get
\begin{align}\label{432}
\widetilde{\Phi}_4&=-i\int_{|\xi'|=1}\int^{+\infty}_{-\infty}{\rm {\rm trace}} [\pi^+_{\xi_n}
\sigma_{0}(\nabla^{\Psi}_{X}\nabla^{\Psi}_{Y}D^{-2}_{\Psi})\times
\partial_{\xi_n}\sigma_{-3}(D^{-2}_{\Psi})](x_0)d\xi_n\sigma(\xi')dx'\nonumber\\
&=i\int_{|\xi'|=1}\int^{+\infty}_{-\infty}{\rm {\rm trace}} [\partial_{\xi_n}\pi^+_{\xi_n}
\sigma_{0}(\nabla^{\Psi}_{X}\nabla^{\Psi}_{Y}D^{-2}_{\Psi})\times
\sigma_{-3}(D^{-2}_{\Psi})](x_0)d\xi_n\sigma(\xi')dx'.
\end{align}
According to Lemma \ref{lem3}, we have
\begin{align}\label{433}
\sigma_{-3}(D^{-2}_{\Psi})(x_0)|_{|\xi'|=1}
=&-\frac{i}{(1+\xi_n^2)^2}\left(-\frac{1}{2}h'(0)\sum_{k<n}\xi_nc(e_k)c(e_n)+\frac{5}{2}h'(0)\xi_n\right)
-\frac{2ih'(0)\xi_n}{(1+\xi_n^2)^3}\nonumber\\
&-\bigg[c(\Psi)ic(\xi)+ic(\xi)c(\Psi)\bigg]|\xi|^{-4}.
\end{align}
In the (\ref{433}), we set $$B_0=-\frac{i}{(1+\xi_n^2)^2}\left(-\frac{1}{2}h'(0)\sum_{k<n}\xi_nc(e_k)c(e_n)+\frac{5}{2}h'(0)\xi_n\right)
-\frac{2ih'(0)\xi_n}{(1+\xi_n^2)^3},$$
and
$$B_1=-\bigg[c(\Psi)ic(\xi)+ic(\xi)c(\Psi)\bigg]|\xi|^{-4},$$
thus we obtain
$$\sigma_{-3}(D^{-2}_{\Psi})(x_0)|_{|\xi'|=1}:=B_0+B_1.$$
By (\ref{432}), we have
\begin{align}\label{434}
&{\rm {\rm trace}} [\partial_{\xi_n}\pi^+_{\xi_n}
\sigma_{0}(\nabla^{\Psi}_{X}\nabla^{\Psi}_{Y}D^{-2}_{\Psi})\times
\sigma_{-3}(D^{-2}_{\Psi})]={\rm {\rm trace}}\bigg[\partial_{\xi_n}\pi^+_{\xi_n}\sigma_{0}(\nabla^{\Psi}_{X}\nabla^{\Psi}_{Y}D^{-2}_{\Psi})\times
\big(B_0+B_1\big)\bigg](x_0).
\end{align}
By Lemma \ref{lem100}, we have
\begin{align}\label{435}
\partial_{\xi_n}\pi^+_{\xi_n}\sigma_{0}(\nabla^{\Psi}_{X}\nabla^{\Psi}_{Y}D^{-2}_{\Psi})(x_0)|_{|\xi'|=1}
=&-\frac{i}{2(\xi_n-i)^2}\sum_{j,l=1}^{n-1}X_jY_l\xi_j\xi_l-\frac{1}{2(\xi_n-i)^2}X_nY_n\nonumber\\
&+\frac{1}{2(\xi_n-i)^2}\sum_{j=1}^{n-1}X_jY_n\xi_j+\frac{1}{2(\xi_n-i)^2}\sum_{l=1}^{n-1}X_nY_l\xi_l.
\end{align}
We observe that $i<n,~\int_{|\xi'|=1}\xi_{i_{1}}\xi_{i_{2}}\cdots\xi_{i_{2d+1}}\sigma(\xi')=0$
and $${\rm {\rm trace}}[c(\Psi)c(\xi)]={\rm {\rm trace}}[c(dx_n)c(\Psi)]\xi_n+\sum\limits^{n-1}_{j=1}{\rm {\rm trace}}[c(dx_j)c(\Psi)]\xi_j',$$
so we omit some items that have no contribution for computing {\bf case (b)}.
Then, we have
\begin{align}\label{436}
&{\rm {\rm trace}}\bigg[\partial_{\xi_n}\pi^+_{\xi_n}\sigma_{0}(\nabla^{\Psi}_{X}\nabla^{\Psi}_{Y}D^{-2}_{\Psi})\times
B_0\bigg](x_0)\nonumber\\
&=-\frac{h'(0)(5\xi_n^2-5+4\xi_n)} {(\xi_n-i)^5(\xi_n+i)^3}\sum_{j,l=1}^{n-1}X_jY_l\xi_j\xi_l
+\frac{h'(0)i(5\xi_n^3-\xi_n)}{(\xi_n-i)^5(\xi_n+i)^3}X_nY_n.
\end{align}
Therefore, we get
\begin{align}\label{437}
&i\int_{|\xi'|=1}\int^{+\infty}_{-\infty}
{\rm {\rm trace}}\bigg[\partial_{\xi_n}\pi^+_{\xi_n}\sigma_{0}(\nabla^{\Psi}_{X}\nabla^{\Psi}_{Y}D^{-2}_{\Psi})\times
B_0\bigg](x_0)d\xi_n\sigma(\xi')dx'\nonumber\\
=&i\int_{|\xi'|=1}\int^{+\infty}_{-\infty}
\bigg(-\frac{h'(0)(5\xi_n^2-5+4\xi_n)} {(\xi_n-i)^5(\xi_n+i)^3}\sum\limits_{j,l=1}^{n-1}X_jY_l\xi_j\xi_l+\frac{h'(0)i(5\xi_n^3-\xi_n)}{(\xi_n-i)^5(\xi_n+i)^3}X_nY_n\bigg)d\xi_n\sigma(\xi')dx'\nonumber\\
=&-i\sum\limits_{j,l=1}^{n-1}X_jY_lh'(0)\Omega_3\int_{\Gamma^{+}}\frac{5\xi_n^2-5+4\xi_n}{(\xi_n-i)^5(\xi_n+i)^2}\xi_j\xi_ld\xi_{n}dx'+iX_nY_nh'(0)\Omega_3\int_{\Gamma^{+}}\frac{5\xi_n^3-\xi_n}{(\xi_n-i)^5(\xi_n+i)^2}d\xi_{n}dx'\nonumber\\
=&-i\sum\limits_{j,l=1}^{n-1}X_jY_lh'(0)\Omega_3\frac{2\pi i}{4!}\left[\frac{5\xi_n^2-5+4\xi_n}{(\xi_n+i)^2}\right]^{(4)}\bigg|_{\xi_n=i}dx'+iX_nY_nh'(0)\Omega_3\frac{2\pi i}{4!}\left[\frac{5\xi_n^3-\xi_n}{(\xi_n+i)^2}\right]^{(4)}\bigg|_{\xi_n=i}dx'\nonumber\\
=&\frac{1}{4}h'(0)\pi\Omega_3dx'\left(\frac{(1-5i)\pi}{3}\sum_{j=1}^{n-1}X_jY_j+\frac{11i}{4}X_nY_n\right).
\end{align}
Similar to (\ref{436}) and (\ref{437}), we have
\begin{align}\label{438}
&i\int_{|\xi'|=1}\int^{+\infty}_{-\infty}
{\rm {\rm trace}}\bigg[\partial_{\xi_n}\pi^+_{\xi_n}\sigma_{0}(\nabla^{\Psi}_{X}\nabla^{\Psi}_{Y}D^{-2}_{\Psi})\times
B_1\bigg](x_0)d\xi_n\sigma(\xi')dx'\nonumber\\
=&\bigg(-i\sum\limits_{j,l=1}^{n-1}X_jY_l-X_nY_n\Omega_3 dx' \bigg)\frac{2\pi
i}{3!}\left[\frac{\xi_n}{2(\xi_n+i)^2}\right]^{(3)}\bigg|_{\xi_n=i}\cdot2{\rm {\rm trace}}[c(\Psi)c(dx_n)]\nonumber\\
=&\left(\frac{\pi^2}{12}\sum_{j=1}^{n-1}X_jY_j-\frac{\pi i}{8}X_nY_n\right)\cdot{\rm {\rm trace}}[c(\Psi)c(dx_n)]\Omega_3 dx'.
\end{align}
From all of this, we get
\begin{align}\label{439}
\widetilde{\Phi}_4
=&\Bigg\{\left(\frac{\pi^2-5i\pi^2}{12}\sum_{j=1}^{n-1}X_jY_j+\frac{11\pi i}{16}X_nY_n\right)h'(0)+\left(\frac{\pi^2}{12}\sum_{j=1}^{n-1}X_jY_j-\frac{\pi i}{8}X_nY_n\right)\nonumber\\
&\times{\rm {\rm trace}}[c(\Psi)c(dx_n)]\Bigg\}\Omega_3 dx'.
\end{align}

\noindent {\bf  case (c)}~$r=-1,~\ell=-2,~k=j=|\alpha|=0$.\\

\noindent By (\ref{fum103}), we get
\begin{align}\label{440}
\widetilde{\Phi}_5=-i\int_{|\xi'|=1}\int^{+\infty}_{-\infty}{\rm {\rm trace}} [\pi^+_{\xi_n}\sigma_{-1}(\nabla_X^{\Psi}\nabla_Y^{\Psi}D^{-2}_{\Psi})\times
\partial_{\xi_n}\sigma_{-2}(D^{-2}_{\Psi})](x_0)d\xi_n\sigma(\xi')dx'.
\end{align}
By Lemma \ref{lem3}, we have
\begin{align}\label{441}
\partial_{\xi_n}\sigma_{-2}(D^{-2}_{\Psi})|_{|\xi'|=1}=-\frac{2\xi_n}{(\xi_n^2+1)^2}.
\end{align}
By Lemma \ref{lem100}, we get
\begin{align}\label{500}
\sigma_{-1}(\nabla^{\Psi}_{X}\nabla^{\Psi}_{Y}D_{\Psi}^{-2})=&
\sigma_{2}(\nabla^{\Psi}_{X}\nabla^{\Psi}_{Y})\sigma_{-3}(D_{\Psi}^{-2})
+\sigma_{1}(\nabla^{\Psi}_{X}\nabla^{\Psi}_{Y})\sigma_{-2}(D_{\Psi}^{-2})\nonumber\\
&+\sum_{j=1}^{n}\partial_{\xi_{j}}\big[\sigma_{2}(\nabla^{\Psi}_{X}\nabla^{\Psi}_{Y})\big]
D_{x_{j}}\big[\sigma_{-2}(D_{\Psi}^{-2})\big].
\end{align}
First, we present the explicit expression for the first term in equation (\ref{500})
\begin{align}
&\sigma_{2}(\nabla^{\Psi}_{X}\nabla^{\Psi}_{Y})\sigma_{-3}(D_{\Psi}^{-2})(x_0)|_{|\xi'|=1}\nonumber\\
=&-\sum_{j,l=1}^{n}X_jY_l\xi_j\xi_l\Bigg[-i|\xi|^{-4}\xi_k(\Gamma^k-2\delta^k)
-i|\xi|^{-6}2\xi^j\xi_\alpha\xi_\beta\partial_jg^{\alpha\beta}-\bigg(c(\Psi)ic(\xi)+ic(\xi)c(\Psi)  \bigg)|\xi|^{-4}\Bigg]\nonumber\\
=&-\sum\limits_{j,l=1}^{n}X_jY_l\xi_j\xi_l\Bigg[
\frac{h'(0)i\sum\limits_{k<n}\xi_nc(e_k)c(e_n)-5i h'(0)\xi_n}{2|\xi|^4}
-\frac{2ih'(0)\xi_n}{(1+\xi_n^2)^3}-\frac{c(\Psi)ic(\xi)+ic(\xi)c(\Psi)}{|\xi|^4}\Bigg].
\end{align}
Next, we present the explicit expression for the second term in equation (\ref{500}),
\begin{align}
&\sigma_{1}(\nabla^{\Psi}_{X}\nabla^{\Psi}_{Y})\sigma_{-2}(D_{\Psi}^{-2})(x_0)|_{|\xi'|=1}\nonumber\\
=&\frac{i\sum\limits_{j,l=1}^nX_j\frac{\partial_{Y_l}}{\partial_{x_j}}i\xi_l
+i\sum\limits_jA(Y)X_j\xi_j+i\sum\limits_lA(Y)Y_l\xi_l+\sum\limits_j  G(X,\Psi)Y_ji\xi_j+\sum\limits_j G(Y,\Psi)X_ji \xi_j}{|\xi|^{2}}.
\end{align}
Finally, we present the explicit expression for the third item of (\ref{500}) in equation (\ref{500}),
\begin{align}
&\sum_{j=1}^{n}\sum_{\alpha}\frac{1}{\alpha!}\partial^{\alpha}_{\xi}\big[\sigma_{2}(\nabla^{\Psi}_{X}\nabla^{\Psi}_{Y})\big]
D_x^{\alpha}\big[\sigma_{-2}(D_{\Psi}^{-2})\big](x_0)|_{|\xi'|=1}
=\sum_{j=1}^{n}\partial_{\xi_{j}}\big[-\sum_{j,l=1}^nX_jY_l\xi_j\xi_l\big]
(-i)\partial_{x_{j}}\big[|\xi|^{-2}\big]\nonumber\\
=&\sum_{j=1}^{n}\sum_{l=1}^{n}i(x_{j}Y_l+x_{l}Y_j)\xi_{l}\partial_{x_{j}}(|\xi|^{-2}).
\end{align}
We  observe that $i<n,~\int_{|\xi'|=1}\xi_{i_{1}}\xi_{i_{2}}\cdots\xi_{i_{2d+1}}\sigma(\xi')=0$ 
, so we omit some items that have no contribution for computing {\bf case (c)}. In addition, we have ${\rm trace}~L(Y)=0$ and ${\rm trace}~[G(X,\Psi)]=-{\rm trace}~[c(X)c(\Psi)]$.
A straightforward computation yields
\begin{align}\label{446}
&{\rm {\rm trace}}[\pi^+_{\xi_n}(\sigma_{2}(\nabla^{\Psi}_{X}\nabla^{\Psi}_{Y})\sigma_{-3}(D_{\Psi}^{-2}))\times
\partial_{\xi_n}\sigma_{-2}(D_{\Psi}^{-2})](x_0)|_{|\xi'|=1}\nonumber\\
=&\frac{h'(0)(2\xi_n^2-\xi_n-2\xi_ni)}{(\xi_n-i)^4(\xi_n+i)^2}\sum_{j,l=1}^{n-1}X_jY_l\xi_j\xi_l\sum_{k<n}\xi_kc(e_k)c(e_n)
+\frac{h'(0)(17\xi_ni-\xi_n^2+4\xi_n^3i)}{(\xi_n-i)^5(\xi_n+i)^2}\sum_{j,l=1}^{n-1}X_jY_l\xi_j\xi_l\nonumber\\
&+\frac{-2i\xi^2_n}{(\xi_n-i)^4(\xi_n+i)^4}\sum_{j,l=1}^{n}X_jY_l\xi_j\xi_l\cdot{\rm 2{\rm trace}}[c(\Psi)c(dx_n)]
\end{align}
Additionally, simple calculations provide
\begin{align}\label{447}
&{\rm {\rm trace}}[\pi^+_{\xi_n}\Big(\sigma_{1}(\nabla^{\Psi}_{X}\nabla^{\Psi}_{Y})\sigma_{-2}(D_{\Psi}^{-2})\Big)\times
\partial_{\xi_n}\sigma_{-2}(D^{-2}_{\Psi})](x_0)|_{|\xi'|=1}\nonumber\\
=&\frac{-\xi_n}{(\xi_n+i)^2(\xi_n-i)^3}\cdot{\rm {\rm trace}}\bigg(iX_n\frac{\partial Y_n}{\partial x_n}+L(Y)X_n+L(X)Y_n+G(X,\Psi)Y_n+G( Y,\Psi)X_n\bigg) ,
\end{align}
and
\begin{align}\label{448}
&-i\int_{|\xi'|=1}\int^{+\infty}_{-\infty}{\rm {\rm trace}}[\pi^+_{\xi_n}\Big(\sigma_{1}(\nabla^{\Psi}_{X}\nabla^{\Psi}_{Y})\sigma_{-2}(D_{\Psi}^{-2})\Big)\times
\partial_{\xi_n}\sigma_{-2}(D^{-2}_{\Psi})](x_0)d\xi_n\sigma(\xi')dx'\nonumber\\
=&\frac{i\pi}{8}\Omega_{3}dx'\cdot{\rm {\rm trace}}\bigg(iX_n\frac{\partial Y_n}{\partial x_n}+L(Y)X_n+L(X)Y_n+G(X,\Psi)Y_n+G( Y,\Psi)X_n\bigg)\nonumber\\
=&\frac{i\pi}{8}\Omega_{3}dx'\cdot\bigg({\rm {\rm trace}}\big(iX_n\frac{\partial Y_n}{\partial x_n}\big)-{\rm trace} [c(X)c(\Psi)]Y_n-{\rm trace} [c(Y)c(\Psi)]X_n\bigg).
\end{align}
In addition, we have
\begin{align}\label{449}
&{\rm {\rm trace}}[\pi^+_{\xi_n}(\sum_{j=1}^{n}\sum_{\alpha}\frac{1}{\alpha!}\partial^{\alpha}_{\xi}
\big[\sigma_{2}(\nabla^{\Psi}_{X}\nabla^{\Psi}_{Y})\big]
D_x^{\alpha}\big[\sigma_{-2}(D_{\Psi}^{-2})\big])\times
\partial_{\xi_n}\sigma_{-2}(D^{-2}_{\Psi})](x_0)|_{|\xi'|=1}\nonumber\\
=&2iX_nY_nh'(0)\xi_n\cdot \frac{-2\xi_n}{(\xi_n^2+1)^2}.
\end{align}

Substituting (\ref{446}),(\ref{447}) and (\ref{449}) into (\ref{440}) yields
\begin{align}\label{450}
\widetilde{\Phi}_5=&\Big(\frac{5i-13}{6}\sum_{j=1}^{n-1}X_jY_j +\frac{3-96i}{8}X_nY_n \Big) h'(0)\pi^2\Omega_3dx'-(\frac{\pi^{2}}{6}\sum_{j=1}^{n-1}X_jY_j\Omega_3dx'
+\frac{\pi}{8}X_nY_n\Omega_3dx')\nonumber\\
&\times{\rm {\rm trace}}[c(\Psi)c(dx_n)]+\frac{i\pi}{8}\Omega_{3}dx'\cdot\bigg({\rm {\rm trace}}\big(iX_n\frac{\partial Y_n}{\partial x_n}\big)-{\rm trace} [c(X)c(\Psi)]Y_n-{\rm trace} [c(Y)c(\Psi)]X_n\bigg).
\end{align}
Let $X=X^T+X_n\partial_n,~Y=Y^T+Y_n\partial_n,$ then we have $\sum\limits_{j=1}^{n-1}X_jY_j=g(X^T,Y^T).$
Now $\widetilde{\Phi}$ is the sum of the cases (a), (b) and (c). Therefore, we get
\begin{align}\label{795}
\widetilde{\Phi}=\sum_{i=1}^5\widetilde{\Phi}_i=&\Big(\frac{15-362i}{32}X_nY_n+\frac{(10i-27)\pi}{24}g(X^T,Y^T)\Big)h'(0)\pi\Omega_3dx'
\nonumber\\
&+\left(\frac{\pi^2}{6}\sum_{j=1}^{n-1}X_jY_j+\frac{\pi-\pi i}{8}X_nY_n\right)\cdot{\rm {\rm trace}}[c(\Psi)c(dx_n)]\Omega_3 dx'\nonumber\\
&+\frac{i\pi}{8}\bigg({\rm {\rm trace}}\big(iX_n\frac{\partial Y_n}{\partial x_n}\big)-{\rm trace} [c(X)c(\Psi)]Y_n-{\rm trace} [c(Y)c(\Psi)]X_n\bigg)\cdot\Omega_{3}dx'
\end{align}
Combining (\ref{fum3.10}) and (\ref{795}), we obtain Theorem \ref{thm1.1}.\\

When $c(\Psi)=f$, we can calculate
${\rm trace}~[c(dx_n)f]=0$ and ${\rm trace}~[c(X)f]=0$,
then we can straightforwardly assert the following corollary:
\begin{cor}
Let $M$ be a $4$-dimensional oriented
compact manifolds with the boundary $\partial M$ and the metric
$g^M$ as above, and let $c(\Psi)=f,$ then
\begin{align*}
&\widetilde{{\rm Wres}}[\pi^{+}\nabla^{\Psi}_{X}\nabla^{\Psi}_{Y}D_{\Psi}^{-2}
    \circ\pi^{+}D_{\Psi}^{-2}]\nonumber\\
=&\frac{4\pi^2}{3}\int_{M}EG(X,Y)vol_{g}+\frac{1}{2}\int_{M}{\rm trace}[\frac{1}{4}s-3f^2]g(X,Y)dvol_{M}\nonumber\\\\
&+\int_{\partial M}\bigg\{\Big(\frac{15-362i}{32}X_nY_n+\frac{(10i-27)\pi}{24}g(X^T,Y^T)\Big)h'(0)\pi\Omega_3
+\frac{i\pi}{8}{\rm {\rm trace}}\bigg(iX_n\frac{\partial Y_n}{\partial x_n}\bigg)\cdot\Omega_{3}\bigg\}d{\rm Vol_{\partial M}},
\end{align*}
where $s$ is the scalar curvature.
\end{cor}

When $c(\Psi)=c(\overline{X})$, we can calculate
\begin{align}\label{4522}
{\rm trace}[c(dx_n)c(\Psi)]=-g(\partial x_n,\overline{X}){\rm trace}[id],
\end{align}
then we can straightforwardly assert the following corollary:
\begin{cor}
Let $M$ be a $4$-dimensional oriented
compact manifolds with the boundary $\partial M$ and the metric
$g^M$ as above, and let $c(\Psi)=c(\overline{X}),$ then
\begin{align}
&\widetilde{{\rm Wres}}[\pi^{+}\nabla^{\Psi}_{X}\nabla^{\Psi}_{Y}D_{\Psi}^{-2}
    \circ\pi^{+}D_{\Psi}^{-2}]=\frac{4\pi^2}{3}\int_{M}EG(X,Y)vol_{g}
+\frac{\pi^2}{2}\int_{M}\Big[-8g(X,\nabla^{TM}_Y\overline{X})
-8g(Y,\nabla^{TM}_X\overline{X})\Big]dvol_{M}\nonumber\\
&+\frac{1}{2}\int_{M}{\rm trace}\big[\frac{1}{4}s+\sum\limits^{n}_{j=1}\frac{1}{2}c(\overline{X})c(e_j)
c(\overline{X})c(e_j)+|\overline{X}|^2\big]g(X,Y)dvol_{M}
+\int_{\partial M}\bigg\{\Big(\frac{15-362i}{32}X_nY_n  \nonumber\\
&+\frac{(10i-27)\pi}{24}g(X^T,Y^T)\Big)h'(0)\pi\Omega_3
-\left(\frac{2\pi^2}{3}\sum_{j=1}^{n-1}X_jY_j+\frac{\pi-\pi i}{2}X_nY_n\right)g(\partial x_n,\overline{X})\Omega_3\nonumber\\
&+\frac{i\pi}{8}{\rm {\rm trace}}\bigg(iX_n\frac{\partial Y_n}{\partial x_n}+4g(X,\overline{X})Y_n+ 4g(Y,\overline{X})X_n\bigg)\cdot\Omega_{3}\bigg\}d{\rm Vol_{\partial M}},
\end{align}
where $s$ is the scalar curvature.
\end{cor}

When $c(\Psi)=c(X_1)c(X_2),$ we have
\begin{align}\label{4555}
{\rm trace}[c(dx_n)c(\Psi)]=0,
\end{align}
and we can compute $\widetilde{{\rm Wres}}[\pi^{+}\nabla^{\Psi}_{X}\nabla^{\Psi}_{Y}D_{\Psi}^{-2}
    \circ\pi^{+}D_{\Psi}^{-2}].$
\begin{cor}
Let $M$ be a $4$-dimensional oriented
compact manifolds with the boundary $\partial M$ and the metric
$g^M$ as above, and let $c(\Psi)=c(X_1)c(X_2),$ then
\begin{align}
&\widetilde{{\rm Wres}}[\pi^{+}\nabla^{\Psi}_{X}\nabla^{\Psi}_{Y}D_{\Psi}^{-2}
    \circ\pi^{+}D_{\Psi}^{-2}]=\frac{4\pi^2}{3}\int_{M}EG(X,Y)vol_{g}\nonumber\\
&+\frac{1}{2}\int_{M}{\rm trace}[\frac{1}{4}s+\sum\limits^{n}_{j=1}\frac{1}{2}c(\Psi)c(e_j)
c(\Psi)c(e_j)-\big(c(\Psi)\big)^2]g(X,Y)dvol_{M}
\nonumber\\
&+\int_{\partial M}\bigg\{\Big(\frac{15-362i}{32}X_nY_n+\frac{(10i-27)\pi}{24}g(X^T,Y^T)\Big)h'(0)\pi\Omega_3
\nonumber\\
&+\frac{i\pi}{8}{\rm {\rm trace}}\big(iX_n\frac{\partial Y_n}{\partial x_n}\big)
\Omega_{3}\bigg\}d{\rm Vol_{\partial M}},
\end{align}
where $s$ is the scalar curvature.
\end{cor}

When $c(\Psi)=c(X_1)c(X_2)c(X_3),$ then we have
\begin{align}\label{4588}
{\rm trace}[c(dx_n)c(\Psi)]=&g(\partial{x_n}, X_1)g(X_2, X_3){\rm trace}[\texttt{id}]-g(\partial{x_n}, X_2)g(X_1, X_3){\rm trace}[\texttt{id}]\nonumber\\
&+g(\partial{x_n}, X_3)g(X_1, X_2){\rm trace}[\texttt{id}],
\end{align}
and we compute that:
\begin{cor}
Let $M$ be a $4$-dimensional oriented
compact manifolds with the boundary $\partial M$ and the metric
$g^M$ as above, and let $c(\Psi)=c(X_1)c(X_2)c(X_3),$ then
\begin{align}
&\widetilde{{\rm Wres}}[\pi^{+}\nabla^{\Psi}_{X}\nabla^{\Psi}_{Y}D_{\Psi}^{-2}
    \circ\pi^{+}D_{\Psi}^{-2}]=\frac{4\pi^2}{3}\int_{M}EG(X,Y)vol_{g}
+\frac{\pi^2}{2}\int_{M}\bigg\{
8\Big[
g(X,X_1)\Big(g(\nabla^{TM}_YX_2,X_3)\nonumber\\&+g(X_2,\nabla^{TM}_YX_3)\Big)
-g(X,X_2)\Big(g(\nabla^{TM}_YX_1,X_3)+g(X_1,\nabla^{TM}_YX_3)\Big)
+g(X,X_3)\Big(g(\nabla^{TM}_YX_1,X_2)\nonumber\\&+g(X_1,\nabla^{TM}_YX_2)\Big)
-g(X,\nabla^{TM}_YX_1)g(X_2,X_3)
-g(X,\nabla^{TM}_YX_2)g(X_1,X_3)
+g(X,\nabla^{TM}_YX_3)g(X_1,X_2)
\Big]\nonumber\\&
-
8\Big[
g(Y,X_1)\Big(g(\nabla^{TM}_XX_2,X_3)+g(X_2,\nabla^{TM}_XX_3)\Big)
-g(Y,X_2)\Big(g(\nabla^{TM}_XX_1,X_3)+g(X_1,\nabla^{TM}_XX_3)\Big)\nonumber\\&
+g(Y,X_3)\Big(g(\nabla^{TM}_XX_1,X_2)+g(X_1,\nabla^{TM}_XX_2)\Big)
-g(Y,\nabla^{TM}_YX_1)g(X_2,X_3)
-g(Y,\nabla^{TM}_YX_2)g(X_1,X_3)\nonumber\\&
+g(Y,\nabla^{TM}_YX_3)g(X_1,X_2)
\Big]
\bigg\}dvol_{M}
+\frac{1}{2}\int_{M}{\rm trace}[\frac{1}{4}s
+\sum\limits^{n}_{j=1}\frac{1}{2}c(\Psi)c(e_j)
c(\Psi)c(e_j)-\big(c(\Psi)\big)^2]g(X,Y)dvol_{M}
\nonumber\\
&+\int_{\partial M}\bigg\{\Big(\frac{15-362i}{32}X_nY_n+\frac{(10i-27)\pi}{24}g(X^T,Y^T)\Big)h'(0)\pi\Omega_3
+\big(\frac{\pi^2}{6}\sum_{j=1}^{n-1}X_jY_j+\frac{\pi-\pi i}{8}X_nY_n\big)\nonumber\\
&\times\bigg(4g(\partial{x_n}, X_1)g(X_2, X_3)-4g(\partial{x_n}, X_2)g(X_1, X_3)+4g(\partial{x_n}, X_3)g(X_1, X_2)\bigg)\Omega_3\nonumber\\
&+\frac{i\pi}{8}\bigg({\rm {\rm trace}}\big(iX_n\frac{\partial Y_n}{\partial x_n}\big)
-4\Big( g(X,X_1)g(X_2,X_3)-g(X,X_2)g(X_1,X_3)+g(X,X_3)g(X_1,X_2) \Big) Y_n\nonumber\\
&
-4\Big( g(Y,X_1)g(X_2,X_3)-g(Y,X_2)g(X_1,X_3)+g(Y,X_3)g(X_1,X_2) \Big)
X_n\bigg)\cdot\Omega_{3}\bigg\}d{\rm Vol_{\partial M}},
\end{align}
where $s$ is the scalar curvature.
\end{cor}

\section{ The noncommutative residue $\widetilde{{\rm Wres}}\bigg[\pi^{+}(\nabla^{\Psi}_{X}\nabla^{\Psi}_{Y}D_{\Psi}^{-1})
    \circ\pi^{+}(D_{\Psi}^{-3})\bigg]$ on manifolds with boundary}
\label{section:4}
In the following, we will compute the residue $\widetilde{{\rm Wres}}\bigg[\pi^{+}(\nabla^{\Psi}_{X}\nabla^{\Psi}_{Y}D_{\Psi}^{-1})
    \circ\pi^{+}(D_{\Psi}^{-3})\bigg]$ on 4 dimensional oriented
compact spin manifolds with boundary and establish a Kastler-Kalau-Walze
type theorem in this case. By (\ref{fum100}), we have
\begin{eqnarray}\label{fum4.1.}
&&\widetilde{{\rm Wres}}\bigg[\pi^{+}(\nabla^{\Psi}_{X}\nabla^{\Psi}_{Y}D_{\Psi}^{-1})
    \circ\pi^{+}(D_{\Psi}^{-3})\bigg]\nonumber\\
&=&\int_{M}\int_{|\xi|=1}{{\rm {\rm trace}}}_{S(TM)}\big[\sigma_{-4}\big( (\nabla^{\Psi}_{X}\nabla^{\Psi}_{Y}D_{\Psi}^{-1})\circ(D_{\Psi}^{-3})\big)\big]\sigma(\xi)dx+\int_{\partial M}\widehat{\Phi},
\end{eqnarray}
where
\begin{eqnarray}\label{fum4.2}
\widehat{\Phi}&=&\int_{|\xi'|=1}\int_{-\infty}^{+\infty}\sum_{j,k=0}^{\infty}\sum \frac{(-i)^{|\alpha|+j+k+\ell}}{\alpha!(j+k+1)!}
{{\rm {\rm trace}}}_{S(TM)}\Big[\partial_{x_{n}}^{j}\partial_{\xi'}^{\alpha}\partial_{\xi_{n}}^{k}\sigma_{r}^{+}
(\nabla^{\Psi}_{X}\nabla^{\Psi}_{Y}D_{\Psi}^{-1})(x',0,\xi',\xi_{n})\nonumber\\
&&\times\partial_{x_{n}}^{\alpha}\partial_{\xi_{n}}^{j+1}\partial_{x_{n}}^{k}\sigma_{l}
(D_{\Psi}^{-3})(x',0,\xi',\xi_{n})\Big]
d\xi_{n}\sigma(\xi')dx' ,
\end{eqnarray}
and the sum is taken over $r-k+|\alpha|+\ell-j-1=-4,r\leq0,\ell\leq-2$.

Locally we can use Theorem \ref{theorem10000} to compute the interior term of (\ref{fum4.1.}), then
 \begin{eqnarray}\label{fum4.3}
&&\int_{M}\int_{|\xi|=1}{{\rm {\rm trace}}}_{S(TM)}\bigg[ \sigma_{-4}\big( (\nabla^{\Psi}_{X}\nabla^{\Psi}_{Y}D_{\Psi}^{-1})\circ(D_{\Psi}^{-3})\big)\bigg]\sigma(\xi)dx\nonumber\\
&=& \frac{4\pi^2}{3}\int_{M}EG(X,Y)vol_{g}
+\frac{\pi^2}{2}\int_{M}{\rm trace}\Big[c(X)\nabla^{S(TM)}_Y\big(c(\Psi)\big)\nonumber\\
&+&\nabla^{S(TM)}_Y\big(c(\Psi)\big) c(X)-c(Y)\nabla^{S(TM)}_X(c(\Psi))-\nabla^{S(TM)}_X(c(\Psi)) c(Y)\Big]dvol_{M}\nonumber\\
&+&\frac{1}{2}\int_{M}{\rm trace}[\frac{1}{4}s+\sum\limits^{n}_{j=1}\frac{1}{2}c(\Psi)c(e_j)
c(\Psi)c(e_j)-\big(c(\Psi)\big)^2]g(X,Y)dvol_{M}
\end{eqnarray}

From Lemma \ref{lem2} and Lemma \ref{lem3}, we have
\begin{lem} \label{lem1000}
The following identities hold:
\begin{align}
\sigma_{1}(\nabla^{\Psi}_{X}\nabla^{\Psi}_{Y}D_{\Psi}^{-1})=&
-i\sum_{j,l=1}^nX_jY_l\xi_j\xi_lc(\xi)|\xi|^{-2};\\
\sigma_{0}(\nabla^{\Psi}_{X}\nabla^{\Psi}_{Y}D_{\Psi}^{-1})=&
\sigma_{2}(\nabla^{\Psi}_{X}\nabla^{\Psi}_{Y})\sigma_{-2}(D_{\Psi}^{-1})
+\sigma_{1}(\nabla^{\Psi}_{X}\nabla^{\Psi}_{Y})\sigma_{-1}(D_{\Psi}^{-1})\nonumber\\
&+\sum_{j=1}^{n}\partial _{\xi_{j}}\big[\sigma_{2}(\nabla^{\Psi}_{X}\nabla^{\Psi}_{Y})\big]
D_{x_{j}}\big[\sigma_{-1}(D_{\Psi}^{-1})\big].
\end{align}
\end{lem}

Write
 \begin{eqnarray}\label{4.2}
D_x^{\alpha}&=(-i)^{|\alpha|}\partial_x^{\alpha};
~\sigma(A^3)=p_3+p_2+p_1+p_0;
~\sigma(A^{-3})=\sum\limits^{\infty}_{j=3}q_{-j}.
\end{eqnarray}
By the composition formula of pseudodifferential operators, we have
\begin{align}\label{4.3}
1=\sigma(A^3\circ A^{-3})
&=\sum_{\alpha}\frac{1}{\alpha!}\partial^{\alpha}_{\xi}[\sigma(A^3)]
A_x^{\alpha}[\sigma(A^{-3})]\nonumber\\
&=(p_3+p_2+p_1+p_0)(q_{-3}+q_{-4}+q_{-5}+\cdots)\nonumber\\
&~~~+\sum_j(\partial_{\xi_j}p_3+\partial_{\xi_j}p_2++\partial_{\xi_j}p_1+\partial_{\xi_j}p_0)
(A_{x_j}q_{-3}+A_{x_j}q_{-4}+A_{x_j}q_{-5}+\cdots)\nonumber\\
&=p_3q_{-3}+(p_3q_{-4}+p_2q_{-3}+\sum_j\partial_{\xi_j}p_3D_{x_j}q_{-3})+\cdots,
\end{align}
so
\begin{equation}\label{4.4}
q_{-3}=p_3^{-1};~q_{-4}=-p_3^{-1}[p_2p_3^{-1}+\sum_j\partial_{\xi_j}p_3D_{x_j}(p_{-3}^{-1})].
\end{equation}
Then it is easy to check that

\begin{lem} \label{lem1001}The following identities hold:
\begin{eqnarray}
\sigma_{-3}(D_{\Psi}^{-3})&=&ic(\xi)|\xi|^{-4};\nonumber\\
\sigma_{-4}(D_{\Psi}^{-3})&=&\frac{c(\xi)\sigma_2(D_{\Psi}^{3})c(\xi)}{|\xi|^8}
+\frac{ic(\xi)}{|\xi|^8}\bigg(|\xi|^4c(dx_n)\partial_{x_n}c(\xi')-2h'(0)c(dx_n)c(\xi)
+2\xi_nc(\xi)\partial{x_n}c(\xi')\nonumber\\
&&+4\xi_nh'(0)\bigg),
\end{eqnarray}
where
 \begin{align}
\sigma_2(D_{\Psi}^{3})=&c(\xi)(4\sigma^k-2\Gamma^k)\xi_{k}
     -\frac{1}{4}|\xi|^2\sum_{s,t}\omega_{s,t}(\widetilde{e_l})c(e_{l})c(\widetilde{e_s})c(\widetilde{e_t})-2|\xi|^2c(\Psi)-2c(\xi)c(\Psi)c(\xi).
\end{align}
\end{lem}

Next, we need to compute $\int_{\partial M}\widehat{\Psi}$.
When $n=4$, then ${\rm trace}_{S(TM)}[{\rm \texttt{id}}]={\rm dim}(\wedge^*(\mathbb{R}^2))=4$, the sum is taken over $
r+l-k-j-|\alpha|=-3,~~r\leq 0,~~l\leq-2,$ then we have the following five cases:\\

\noindent {\bf case (1)}~$r=1,~l=-2,~k=j=0,~|\alpha|=1$.\\

\noindent By (\ref{fum4.2}), we get
\begin{equation}
\label{b24}
\widehat{\Phi}_1=-\int_{|\xi'|=1}\int^{+\infty}_{-\infty}\sum_{|\alpha|=1}
{\rm {\rm trace}}[\partial^\alpha_{\xi'}\pi^+_{\xi_n}\sigma_{1}(\nabla^{\Psi}_{X}\nabla^{\Psi}_{Y}D_{\Psi}^{-1})\times
 \partial^\alpha_{x'}\partial_{\xi_n}\sigma_{-3}(D_{\Psi}^{-3})](x_0)d\xi_n\sigma(\xi')dx'.
\end{equation}
By Lemma 2.2 in \cite{Wa3}, for $i<n$, then
\begin{equation}
\label{b25}
\partial_{x_i}\sigma_{-3}(D_{\Psi}^{-3})(x_0)=
\partial_{x_i}(ic(\xi)|\xi|^{-4})(x_0)=
i\frac{\partial_{x_i}c(\xi)}{|\xi|^4}(x_0)
+i\frac{c(\xi)\partial_{x_i}(|\xi|^4)}{|\xi|^8}(x_0)
=0,
\end{equation}
\noindent so $\widehat{\Phi}_1=0$.\\

\noindent {\bf case (2)}~$r=1,~l=-3,~k=|\alpha|=0,~j=1$.\\

\noindent  By (\ref{fum4.2}), we get
\begin{equation}
\label{b26}
\widehat{\Phi}_2=-\frac{1}{2}\int_{|\xi'|=1}\int^{+\infty}_{-\infty}
{\rm {\rm trace}} [\partial_{x_n}\pi^+_{\xi_n}\sigma_{1}(\nabla^{\Psi}_{X}\nabla^{\Psi}_{Y}D_{\Psi}^{-1})\times
\partial_{\xi_n}^2\sigma_{-3}(D_{\Psi}^{-3})](x_0)d\xi_n\sigma(\xi')dx'.
\end{equation}
By Lemma \ref{lem1001}, we have
\begin{eqnarray}\label{b237}
\partial_{\xi_n}^2\sigma_{-3}(D_{\Psi}^{-3})(x_0)=\partial_{\xi_n}^2(c(\xi)|\xi|^{-4})(x_0)
=\frac{(20\xi_n^2-4)ic(\xi')+12(\xi^3-\xi)ic(dx_n)}{(1+\xi_n^2)^4},
\end{eqnarray}
In addition, by Lemma \ref{lem1000}, we get
\begin{align}\label{b27}
&\partial_{x_n}\sigma_{1}(\nabla^{\Psi}_{X}\nabla^{\Psi}_{Y}D_{\Psi}^{-1})(x_0)=\sum_{j,l=1}^nX_jY_l\xi_j\xi_l \left[ \frac{\partial_{x_n}c(\xi')}{1+\xi_n^2}+\frac{c(\xi)h'(0)|\xi'|^2}{(1+\xi_n^2)^2}\right].
\end{align}
Then, we have
\begin{align}\label{b28}
&\pi^+_{\xi_n}\partial_{x_n}\sigma_{1}(\nabla^{\Psi}_{X}\nabla^{\Psi}_{Y}D_{\Psi}^{-1})(x_0)
=\partial_{x_n}\pi^+_{\xi_n}\sigma_{1}(\nabla^{\Psi}_{X}\nabla^{\Psi}_{Y}D_{\Psi}^{-1})(x_0)\nonumber\\
=&i\sum_{j,l=1}^{n-1}X_jY_l\xi_j\xi_lh'(0)|\xi'|^2\left[\frac{ic(\xi')}{4(\xi_n-i)}+\frac{c(\xi')
+ic(dx_n)}{4(\xi_n-i)^2}\right]-\sum_{j,l=1}^{n-1}X_jY_l\xi_j\xi_l\frac{\partial_{x_n}c(\xi')}{2(\xi_n-i)}
-iX_nY_n\left\{\frac{\partial_{x_n}c(\xi')}{2(\xi_n-i)}\right.\nonumber\\
&\left.+h'(0)|\xi'|\left[-\frac{2ic(\xi')-3c(dx_n)}{4(\xi_n-i)}
\frac{[c(\xi')+ic(dx_n)][i(\xi_n-i)+1]}{4(\xi_n-i)^2}\right]\right\}
-\sum_{j=1}^{n-1}X_jY_n\xi_j\left[i\frac{\partial_{x_n}c(\xi')}{2(\xi_n-i)}\right.\nonumber\\
&\left.
-\frac{ih'(0)|\xi'|[c(\xi')+2ic(dx_n)]}{4(\xi_n-i)}
-\frac{[ic(\xi')-c(dx_n)][i(\xi_n-i)+1]}{(\xi_n-i)^2}\right]-\sum_{l=1}^{n-1}X_nY_l\xi_l \left[i\frac{\partial_{x_n}c(\xi')}{2(\xi_n-i)}\right.\nonumber\\
&\left.-\frac{ih'(0)|\xi'|[c(\xi')+2ic(dx_n)]}{4(\xi_n-i)}
-\frac{[ic(\xi')-c(dx_n)][i(\xi_n-i)+1]}{(\xi_n-i)^2}\right] .
\end{align}
We note that $i<n,~\int_{|\xi'|=1}\xi_{i_{1}}\xi_{i_{2}}\cdots\xi_{i_{2r+1}}\sigma(\xi')=0$,
so we omit some items that have no contribution for computing {\bf case (2)}.
Then there is the following formula
\begin{align}\label{33}
&{\rm {\rm trace}}[\partial_{x_n}\pi^+_{\xi_n}\sigma_{1}(\nabla^{\Psi}_{X}\nabla^{\Psi}_{Y}D_{\Psi}^{-1})\times
\partial_{\xi_n}^2\sigma_{-3}(D_{\Psi}^{-3})](x_0)\nonumber\\
=&\sum\limits_{j,l=1}^{n-1}X_jY_l\xi_j\xi_lh'(0)\left[8i\frac{5\xi_n^2-1}{(\xi_n-i)^5(\xi_n+i)^4}
 +\frac{4(5\xi_n^2-1)+12i(\xi_n^3-\xi_n)}{(\xi_n-i)^6(\xi_n+i)^4}\right]\nonumber\\
&+X_nY_nh'(0)\left[\frac{(4i-4)(\xi^2-1)+48(\xi_n^3-\xi_n)}{(\xi_n-i)^5(\xi_n+i)^4}
 -\frac{4(5\xi_n^2-1)+12i(\xi_n^3-\xi_n)}{(\xi_n-i)^6(\xi_n+i)^4}\right]\nonumber\\
&+8\sum\limits_{j=1}^{n-1}X_jY_n\xi_j\left[\frac{(6-3ih'(0))(\xi_n^3-\xi_n)-2i(5\xi_n^2-1)}{(\xi_n-i)^5(\xi_n+i)^4}
 +\frac{2(5\xi_n^2-1)+6i(\xi_n^3-\xi_n)}{(\xi_n-i)^6(\xi_n+i)^4}\right]\nonumber\\
&+8\sum\limits_{l=1}^{n-1}X_nY_l\xi_l\left[\frac{(6-3ih'(0))(\xi_n^3-\xi_n)-2i(5\xi_n^2-1)}{(\xi_n-i)^5(\xi_n+i)^4}
 +\frac{2(5\xi_n^2-1)+6i(\xi_n^3-\xi_n)}{(\xi_n-i)^6(\xi_n+i)^4}\right].
\end{align}
Therefore, we get
\begin{align}\label{35}
\widehat{\Phi}_2=&\frac{1}{2}\int_{|\xi'|=1}\int^{+\infty}_{-\infty}\bigg\{
\sum_{j,l=1}^{n-1}X_jY_l\xi_j\xi_lh'(0)\left[8i\frac{5\xi_n^2-1}{(\xi_n-i)^5(\xi_n+i)^4}
 +\frac{4(5\xi_n^2-1)+12i(\xi_n^3-\xi_n)}{(\xi_n-i)^6(\xi_n+i)^4}\right]\nonumber\\
&+X_nY_nh'(0)\left[\frac{(4i-4)(\xi^2-1)+48(\xi_n^3-\xi_n)}{(\xi_n-i)^5(\xi_n+i)^4}
 -\frac{4(5\xi_n^2-1)+12i(\xi_n^3-\xi_n)}{(\xi_n-i)^6(\xi_n+i)^4}\right]
 \bigg\}d\xi_n\sigma(\xi')dx'\nonumber\\
=&\sum_{j,l=1}^{n-1}X_jY_lh'(0)\Omega_3\int_{\Gamma^{+}}\left[8i\frac{5\xi_n^2-1}{(\xi_n-i)^5(\xi_n+i)^4}
 +\frac{4(5\xi_n^2-1)+12i(\xi_n^3-\xi_n)}{(\xi_n-i)^6(\xi_n+i)^4}\right]\xi_j\xi_ld\xi_{n}dx'\nonumber\\
 &+X_nY_nh'(0)\Omega_3\int_{\Gamma^{+}}\left[\frac{(4i-4)(\xi^2-1)+48(\xi_n^3-\xi_n)}{(\xi_n-i)^5(\xi_n+i)^4}
 -\frac{4(5\xi_n^2-1)+12i(\xi_n^3-\xi_n)}{(\xi_n-i)^6(\xi_n+i)^4}\right]d\xi_{n}dx'\nonumber\\
=&-\left[\frac{592}{3}\pi\sum_{j=1}^{n-1}X_jY_j
+\left(\frac{461}{4}+\frac{23}{4}i\right)X_nY_n\right]h'(0)\pi\Omega_3dx',
\end{align}
where ${\rm \Omega_{3}}$ is the canonical volume of $S^{2}.$\\

\noindent {\bf case (3)}~$r=1,~l=-3,~j=|\alpha|=0,~k=1$.\\

 \noindent By (\ref{fum4.2}), we get
\begin{align}\label{36}
\widehat{\Phi}_3&=-\frac{1}{2}\int_{|\xi'|=1}\int^{+\infty}_{-\infty}
{\rm {\rm trace}} [\partial_{\xi_n}\pi^+_{\xi_n}\sigma_{1}(\nabla^{\Psi}_{X}\nabla^{\Psi}_{Y}D_{\Psi}^{-1})\times
\partial_{\xi_n}\partial_{x_n}\sigma_{-3}(D_{\Psi}^{-3})](x_0)d\xi_n\sigma(\xi')dx'\nonumber\\
&=\frac{1}{2}\int_{|\xi'|=1}\int^{+\infty}_{-\infty}
{\rm {\rm trace}} [\partial_{\xi_n}^2\pi^+_{\xi_n}\sigma_{1}(\nabla^{\Psi}_{X}\nabla^{\Psi}_{Y}D_{\Psi}^{-1})\times
\partial_{x_n}\sigma_{-3}(D_{\Psi}^{-3})](x_0)d\xi_n\sigma(\xi')dx'.
\end{align}
\noindent By Lemma 5.3, we have
\begin{eqnarray}\label{37}
\partial_{x_n}\sigma_{-3}(D_{\Psi}^{-3})(x_0)|_{|\xi'|=1}
=\frac{i\partial_{x_n}[c(\xi')]}{(1+\xi_n^2)^4}-\frac{2ih'(0)c(\xi)|\xi'|^2_{g^{\partial M}}}{(1+\xi_n^2)^6}.
\end{eqnarray}
By integrating formula we obtain
\begin{align}\label{38}
\pi^+_{\xi_n}\sigma_{1}(\nabla^{\Psi}_{X}\nabla^{\Psi}_{Y}D_{\Psi}^{-1})
=&-\frac{c(\xi')+ic(dx_n)}{2(\xi_n-i)}\sum\limits_{j,l=1}^{n-1}X_jY_l\xi_j\xi_l
-\frac{c(\xi')+ic(dx_n)}{2(\xi_n-i)}X_nY_n\nonumber\\
&-\frac{ic(\xi')-c(dx_n)}{2(\xi_n-i)}\sum\limits_{j=1}^{n-1}X_jY_n\xi_j
-\frac{ic(\xi')-c(dx_n)}{2(\xi_n-i)}\sum\limits_{l=1}^{n-1}X_nY_l\xi_l.
\end{align}
Then, we have
\begin{align}\label{mmmmm}
\partial_{\xi_n}^2\pi^+_{\xi_n}\sigma_{1}(\nabla^{\Psi}_{X}\nabla^{\Psi}_{Y}D_{\Psi}^{-1})
=-\frac{c(\xi')+ic(dx_n)}{(\xi_n-i)^3}\sum\limits_{j,l=1}^{n-1}X_jY_l\xi_j\xi_l
-\frac{c(\xi')+ic(dx_n)}{(\xi_n-i)^3}X_nY_n.
\end{align}

We note that $i<n,~\int_{|\xi'|=1}\xi_{i_{1}}\xi_{i_{2}}\cdots\xi_{i_{2r+1}}\sigma(\xi')=0$,
so we omit some items that have no contribution for computing {\bf case (3)}, then
\begin{align}\label{39}
&{\rm {\rm trace}} [\partial_{\xi_n}\pi^+_{\xi_n}\sigma_{1}(\nabla^{\Psi}_{X}\nabla^{\Psi}_{Y}D_{\Psi}^{-1})\times
\partial_{\xi_n}\partial_{x_n}\sigma_{-3}(D_{\Psi}^{-3})](x_0)=\frac{-2 h'(0)}{(\xi_n-i)^5(\xi_n+i)^2}(\sum\limits_{j,l=1}^{n-1}X_jY_l\xi_j\xi_l +X_nY_n).
\end{align}
Therefore, we get
\begin{align}\label{41}
\widehat{\Phi}_3&=\frac{1}{2}\int_{|\xi'|=1}\int^{+\infty}_{-\infty}
\bigg(\frac{-2 h'(0)}{(\xi_n-i)^5(\xi_n+i)^2}(\sum\limits_{j,l=1}^{n-1}X_jY_l\xi_j\xi_l +X_nY_n)\bigg)d\xi_n\sigma(\xi')dx'\nonumber\\
&=h'(0)\Omega_3\bigg(-\sum_{j,l=1}^{n-1}X_jY_l\int_{\Gamma^{+}}\frac{1}{(\xi_n-i)^5(\xi_n+i)^2}\xi_j\xi_ld\xi_{n}dx'
-X_nY_n\int_{\Gamma^{+}}\frac{1}{(\xi_n-i)^5(\xi_n+i)^2}d\xi_{n}dx'\bigg)\nonumber\\
&=\bigg(-\sum_{j,l=1}^{n-1}X_jY_l+2X_nY_n\bigg)h'(0)\Omega_3\frac{2\pi i}{4!}
\left[\frac{1}{(\xi_n+i)^2}\right]^{(4)}\bigg|_{\xi_n=i}dx'\nonumber\\
&=\left(\frac{5\pi i}{6}\sum\limits_{j=1}^{n-1}X_jY_j+\frac{5i}{8}X_nY_n\right)h'(0)\pi\Omega_3dx'.
\end{align}

\noindent {\bf case (4)}~$r=0,~l=-3,~k=j=|\alpha|=0$.\\

\noindent
By (\ref{fum4.2}), we get
\begin{align}\label{522}
\widehat{\Phi}_4&=-i\int_{|\xi'|=1}\int^{+\infty}_{-\infty}
{\rm {\rm trace}}[\pi^+_{\xi_n}\sigma_{0}(\nabla^{\Psi}_{X}\nabla^{\Psi}_{Y}D_{\Psi}^{-1})\times
\partial_{\xi_n}\sigma_{-3}(D_{\Psi}^{-3})](x_0)d\xi_n\sigma(\xi')dx'.
\end{align}
By Lemma 5.3, we obtain
\begin{align}\label{523}
\partial_{\xi_n}\sigma_{-3}(D_{\Psi}^{-3})(x_0)|_{|\xi'|=1}
=\frac{ic(dx_n)}{(1+\xi_n^2)^2}-\frac{4i\xi_nc(\xi)}{(1+\xi_n^2)^3}.
\end{align}
By Lemma 5.2, we have
\begin{align}\label{524}
\sigma_{0}(\nabla^{\Psi}_{X}\nabla^{\Psi}_{Y}D_{\Psi}^{-1})=&
\sigma_{2}(\nabla^{\Psi}_{X}\nabla^{\Psi}_{Y})\sigma_{-2}(D_{\Psi}^{-1})
+\sigma_{1}(\nabla^{\Psi}_{X}\nabla^{\Psi}_{Y})\sigma_{-1}(D_{\Psi}^{-1})\nonumber\\
&+\sum_{j=1}^{n}\partial _{\xi_{j}}\big[\sigma_{2}(\nabla^{\Psi}_{X}\nabla^{\Psi}_{Y})\big]
D_{x_{j}}\big[\sigma_{-1}(D_{\Psi}^{-1})\big].
\end{align}
(1)First, we present the explicit expression for the first term in equation (\ref{524})
\begin{align}
&\sigma_{2}(\nabla^{\Psi}_{X}\nabla^{\Psi}_{Y})\sigma_{-2}(D_{\Psi}^{-1})(x_0)|_{|\xi'|=1}\nonumber\\
=&-\sum_{j,l=1}^{n}X_jY_l\xi_j\xi_l\left[\frac{c(\xi)\sigma_0(D_{\Psi})c(\xi)}{|\xi|^4}
+\frac{c(\xi)}{|\xi|^6}\sum\limits_jc(dx_j)\Big(\partial_{x_j}[c(\xi)]|\xi|^2-c(\xi)\partial_{x_j}(|\xi|^2)\Big)\right],
\end{align}
By integrating formula, we obtain
\begin{align}\label{526}
&\pi^+_{\xi_n}\bigg\{-\sum\limits_{j,l=1}^{n}X_jY_l\xi_j\xi_l
\left[\frac{c(\xi)\sigma_0(D_{\Psi})(x_0)c(\xi)+c(\xi)c(dx_n)\partial_{x_n}[c(\xi')](x_0)}{(1+\xi_n^2)^2}\right]
\bigg\}
\nonumber\\
=&-\frac{\sum\limits_{j,l=1}^{n-1}X_jY_l\xi_j\xi_l\cdot K_1}{4(\xi_n-i)}-\frac{\sum\limits_{j,l=1}^{n-1}X_jY_l\xi_j\xi_l\cdot K_2}{4(\xi_n-i)^2}+\frac{K_3}{4(\xi_n-i)^2},
\end{align}
where
\begin{align}
K_1=&ic(\xi')\sigma_0(D)c(\xi')+ic(dx_n)(-\frac{3}{4}h'(0)c(dx_n))c(dx_n)+ic(\xi')c(dx_n)\partial_{x_n}[c(\xi')],\\
K_2=&[c(\xi')+ic(dx_n)]\sigma_0(D)[c(\xi')+ic(dx_n)]+c(\xi')c(dx_n)\partial_{x_n}c(\xi')-i\partial_{x_n}[c(\xi')],\\
K_3=&\sum_{j,l=1}^{n-1}X_{j}Y_l\xi_{j}\xi_{l}\Big( (-2-i\xi_{n})c(\xi')c(\Psi)c(\xi')-ic(dx_n)c(\Psi)c(\xi')-i c(\xi')c(\Psi)c(dx_n)\nonumber\\
   &-i\xi_{n}c(dx_n)c(\Psi)c(dx_n)\Big)+X_{n}Y_n\Big( -i\xi_{n}c(\xi')c(\Psi)c(\xi')-ic(dx_n)c(\Psi)c(\xi')\nonumber\\
   &-i c(\xi')c(\Psi)c(dx_n)+c(dx_n)c(\Psi)c(dx_n)\Big).
\end{align}
Using a similar approach, we derive
\begin{align}
\pi^+_{\xi_n}\left[\frac{c(\xi)c(dx_n)c(\xi)}{(1+\xi_n)^3}(x_0)|_{|\xi'|=1}\right]
&=\frac{1}{2}\left[\frac{c(dx_n)}{4i(\xi_n-i)}+\frac{c(dx_n)-ic(\xi')}{8(\xi_n-i)^2}
+\frac{3\xi_n-7i}{8(\xi_n-i)^3}[ic(\xi')-c(dx_n)]\right].
\end{align}
(2)Next, we present the explicit expression for the second term in equation (5.28)
\begin{align}
&\sigma_{1}(\nabla^{\Psi}_{X}\nabla^{\Psi}_{Y})\sigma_{-1}(D_{\Psi}^{-1})(x_0)|_{|\xi'|=1}\nonumber\\
=&\Big(i\sum_{j,l=1}^nX_j\frac{\partial_{Y_l}}{\partial_{x_j}}\partial_{x_l}
+i\sum_jL(Y)X_j\xi_j+i\sum_lL(Y)Y_l\xi_l+\sum_ji G(X,\Psi)Y_j\xi_j\nonumber\\
&+\sum_ji G(Y,\Psi)X_j\xi_j \Big)\frac{ic(\xi)}{|\xi|^{2}};
\end{align}
By integrating formula we get
\begin{align}
&\pi^+_{\xi_n}\left(\Big(\sum_j^{n}i G(X,\Psi)Y_j\xi_j+\sum_j^{n}i G(Y,\Psi)X_j\xi_j \Big)
\frac{ic(\xi)}{|\xi|^{2}}\right)\nonumber\\
=&\pi^+_{\xi_n}\left(\Big(\sum_j^{n-1}i G(X,\Psi)Y_j\xi_j+\sum_j^{n-1}i G(Y,\Psi)X_j\xi_j \Big)
\frac{ic(\xi)}{|\xi|^{2}}\right)\nonumber\\
&+\pi^+_{\xi_n}\left(\Big( i G(X,\Psi)Y_n\xi_n+i G(Y,\Psi)X_n\xi_n \Big)
\frac{ic(\xi)}{|\xi|^{2}}\right)\nonumber\\
=&\Big(\sum_j^{n-1}i G(X,\Psi)Y_j\xi_j+\sum_j^{n-1}i G(Y,\Psi)X_j\xi_j \Big)
\frac{ic(\xi')-c(dx_n)}{2(\xi_{n}-i)}\nonumber\\
&+\Big( i G(X,\Psi)Y_n\xi_n+iG(Y,\Psi)X_n\xi_n \Big)
\frac{-c(\xi')-ic(dx_n)}{2(\xi_{n}-i)}
\end{align}
We note that $i<n,~\int_{|\xi'|=1}\xi_{i_{1}}\xi_{i_{2}}\cdots\xi_{i_{2d+1}}\sigma(\xi')=0$,
and ${\rm trace}[c(\xi')]={\rm trace}[c(dx_{n})]=-4$,
then
\begin{align}\label{39}
&-i\int_{|\xi'|=1}\int^{+\infty}_{-\infty}{\rm trace} \left(\pi^+_{\xi_n}\Big(\sigma_{1}(\nabla^{\Psi}_{X}\nabla^{\Psi}_{Y})\sigma_{-1}(D_{\Psi}^{-1})\Big)\times
\partial_{\xi_n}\sigma_{-3}(D_{\Psi}^{-3})\right)(x_0)(x_0)d\xi_n\sigma(\xi')dx'\nonumber\\
=&\frac{3\pi}{4}\Big[i X_n \cdot\frac{\partial Y_n}{\partial x_n}+L(Y)X_n+L(X)Y_n+{\rm trace}[G(X,\Psi)]Y_n+{\rm trace}[G(Y,\Psi)]X_n\Big].
\end{align}
(3)
Finally, we present the explicit expression for the second term in equation (5.28)
\begin{align}
&\sum_{j=1}^{n}\sum_{\alpha}\frac{1}{\alpha!}\partial^{\alpha}_{\xi}\big[\sigma_{2}(\nabla^{\Psi}_{X}\nabla^{\Psi}_{Y})\big]
D_{x_{j}}\big[\sigma_{-1}(D_{\Psi}^{-1})\big](x_0)|_{|\xi'|=1}\nonumber\\
=&\sum_{j=1}^{n}\partial_{\xi_{j}}
\big[\sigma_{2}(\nabla^{\Psi}_{X}\nabla^{\Psi}_{Y})\big]
(-i)\partial_{x_{j}}\big[\sigma_{-1}(D_{\Psi}^{-1})\big]\nonumber\\
=&\sum_{j=1}^{n}\partial_{\xi_{j}}\big[-\sum_{j,l=1}^nX_jY_l\xi_j\xi_l\big]
(-i)\partial_{x_{j}}\big[\frac{ic(\xi)}{|\xi|^{2}}\big]\nonumber\\
=&\sum_{j=1}^{n}\sum_{l=1}^{n}i(x_{j}Y_l+x_{l}Y_j)\xi_{l}\partial_{x_{j}}(\frac{ic(\xi)}{|\xi|^{2}}).
\end{align}
By integrating formula we obtain
\begin{align}\label{547}
&\pi^+_{\xi_n}\left(\sum_{j=1}^{n}\sum_{\alpha}\frac{1}{\alpha!}\partial^{\alpha}_{\xi}\big[\sigma_{2}(\nabla^{\Psi}_{X}\nabla^{\Psi}_{Y})\big]
D_{x_{j}}\big[\sigma_{-1}(D_{\Psi}^{-1})\big]\right)\nonumber\\
=&\pi^+_{\xi_n}\left(\sum_{l=1}^{n-1}i(x_{n}Y_l+x_{l}Y_n)\xi_{l}\partial_{x_{n}}(\frac{ic(\xi)}{|\xi|^{2}})
\right)+\pi^+_{\xi_n}\left(i(x_{n}Y_n+x_{n}Y_n)\xi_{n}\partial_{x_{n}}(\frac{ic(\xi)}{|\xi|^{2}})
\right)\nonumber\\
=&\sum_{l=1}^{n-1}(x_{n}Y_l+x_{l}Y_n)\xi_{l}\Big(\frac{i\partial_{x_{n}}(c(\xi'))}{2(\xi_n-i)}
+h'(0)\frac{(-2-i\xi_n )c(\xi')}{4(\xi_n-i)^2}
-h'(0)\frac{ic(dx_n)}{4(\xi_n-i)^2}\Big)\nonumber\\
&+x_{n}Y_n\Big(\frac{-\partial_{x_{n}}(c(\xi'))}{(\xi_n-i)}
+h'(0)\frac{(-i )c(\xi')}{2(\xi_n-i)^2}
-h'(0)\frac{-i\xi_n c(dx_n)}{2(\xi_n-i)^2}\Big)
\end{align}
Substituting (\ref{523}) and (\ref{547}) into (\ref{522}) yields
\begin{align}\label{39}
&-i\int_{|\xi'|=1}\int^{+\infty}_{-\infty}
{\rm trace}\Big[\pi^+_{\xi_n}\left(\sum_{j=1}^{n}\sum_{\alpha}\frac{1}{\alpha!}\partial^{\alpha}_{\xi}\big[\sigma_{2}(\nabla^{\Psi}_{X}\nabla^{\Psi}_{Y})\big]
D_{x_{j}}\big[\sigma_{-1}(D_{\Psi}^{-1})\big]\right)\partial_{\xi_n}\sigma_{-3}(D_{\Psi}^{-3})\Big](x_0)\nonumber\\
&\times
d\xi_n\sigma(\xi')dx'\nonumber\\
=&\frac{7-15i}{8}X_{n}Y_n\pi h'(0)\Omega_3dx'.
\end{align}
Combining (1), (2), and (3) yields the desired equality
\begin{align}\label{41}
\widehat{\Phi}_4
=&\Big(\frac{55\pi}{24}\sum_{j=1}^{n-1}X_jY_j
+\frac{15i-60+\pi}{8}X_nY_n\Big)h'(0)\pi\Omega_3dx'
-\Big(\frac{\pi^2}{3}\sum_{j,l=1}^{n-1}X_{j}Y_l
+\frac{3\pi}{4}X_{n}Y_n\Big)
{\rm trace}[c(dx_n)c(\Psi)]\Omega_3dx'\nonumber\\
&-\frac{3\pi}{4}X_{n}Y_n{\rm trace}[c(dx_n)\sigma_0(D)]\Omega_3dx'+\frac{3\pi i}{16}{\rm trace}[-c(X)c(\Psi)X_n-c(Y)c(\Psi)Y_n]\Omega_3dx'
.
\end{align}

\noindent {\bf  case (5)}~$r=1,~\ell=-4,~k=j=|\alpha|=0$.\\

\noindent
By (\ref{fum4.2}), we get
\begin{align}\label{61}
\widehat{\Phi}_5&=-\int_{|\xi'|=1}\int^{+\infty}_{-\infty}{\rm trace} [\pi^+_{\xi_n}
\sigma_{1}(\nabla^{\Psi}_{X}\nabla^{\Psi}_{Y}D_{\Psi}^{-1})\times
\partial_{\xi_n}\sigma_{-4}(D_{\Psi}^{-3})](x_0)d\xi_n\sigma(\xi')dx'\nonumber\\
&=\int_{|\xi'|=1}\int^{+\infty}_{-\infty}{\rm trace}
[\partial_{\xi_n}\pi^+_{\xi_n}\sigma_{1}(\nabla^{\Psi}_{X}\nabla^{\Psi}_{Y}D_{\Psi}^{-1})\times
\sigma_{-4}(D_{\Psi}^{-3})](x_0)d\xi_n\sigma(\xi')dx'.
\end{align}
By Lemma \ref{lem1001}, we have
\begin{align}\label{62}
\sigma_{-4}(D_{\Psi}^{-3})(x_0)|_{|\xi'|=1}=&\frac{1}{(\xi_n^2+1)^4}
\left[\left(\frac{11}{2}\xi_n(1+\xi_n^2)+8i\xi_n\right)h'(0)c(\xi')\right.\nonumber\\
&+\left[-2i+6i\xi_n^2-\frac{7}{4}(1+\xi_n^2)
 +\frac{15}{4}\xi_n^2(1+\xi^2_n)\right]h'(0)c(dx_n) \nonumber\\
&-3i\xi_n(1+\xi^2_n)\partial_{x_n}c(\xi')
 +i(1+\xi^2_n)c(\xi')c(dx_n)\partial_{x_n}c(\xi') \nonumber\\
 &+2\Big[\frac{c(\xi)c(\Psi)c(\xi)}{|\xi|^6}-\frac{c(\Psi)}{|\xi|^4}\Big],
\end{align}
and
\begin{align}\label{621}
\partial_{\xi_n}\pi^+_{\xi_n}\sigma_{1}(\nabla^{\Psi}_{X}\nabla^{\Psi}_{Y}D_{\Psi}^{-1})(x_0)|_{|\xi'|=1}
&=\frac{c(\xi')+ic(dx_n)}{2(\xi_n-i)^2}\sum\limits_{j,l=1}^{n-1}X_jY_l\xi_j\xi_l
-\frac{c(\xi')+ic(dx_n)}{2(\xi_n-i)^2}X_nY_n \nonumber\\
&+\frac{ic(\xi')-c(dx_n)}{2(\xi_n-i)^2}\sum\limits_{j=1}^{n}X_jY_n\xi_j
+\frac{ic(\xi')-c(dx_n)}{2(\xi_n-i)^2}\sum\limits_{l=1}^{n}X_nY_l\xi_l.
\end{align}

We note that $i<n,~\int_{|\xi'|=1}\xi_{i_{1}}\xi_{i_{2}}\cdots\xi_{i_{2d+1}}\sigma(\xi')=0$,
so we omit some items that have no contribution for computing {\bf case (5)}. Here
\begin{align}\label{63}
{\rm trace}[c(\xi')c(\xi')c(dx_n)\partial_{x_n}c(\xi')]=0 ;~
{\rm trace}[c(dx_n)c(\xi')c(dx_n)\partial_{x_n}c(\xi')]=-2h'(0).
\end{align}
Also, straightforward computations yield
\begin{align}\label{71}
&{\rm trace}\bigg[\partial_{\xi_n}\pi^+_{\xi_n}\sigma_{-1}(\nabla^{\Psi}_{X}\nabla^{\Psi}_{Y}D_{\Psi}^{-1})\times
\frac{1}{(\xi_n^2+1)^4}
\Big(\big(\frac{11}{2}\xi_n(1+\xi_n^2)+8i\xi_n\big)h'(0)c(\xi')\nonumber\\
&+\big(-2i+6i\xi_n^2-\frac{7}{4}(1+\xi_n^2)
 +\frac{15}{4}\xi_n^2(1+\xi^2_n)\big)h'(0)c(dx_n) \nonumber\\
&-3i\xi_n(1+\xi^2_n)\partial_{x_n}c(\xi')
 +i(1+\xi^2_n)c(\xi')c(dx_n)\partial_{x_n}c(\xi')\Big) \bigg](x_0)|_{|\xi'|=1}\nonumber\\
=&\sum_{j,l=1}^{n-1}X_jY_l\xi_j\xi_l\frac{h'(0)(7+6i-(20-15i)\xi_n-(7-6i)\xi_n^2+15i\xi_n^3)}{(\xi_n-i)^5(\xi_n+i)^4}\nonumber\\
&+X_nY_n\frac{(3i-11)\xi_n(1-\xi_n^2)-16i\xi_n+(13+\frac{7}{2}i)(1+\xi_n^2)-16-\frac{15}{2}\xi_n^2(1+\xi_n^2)}{(\xi_n-i)^2(\xi_n+i)^4},
\end{align}
and
\begin{align}\label{710}
&{\rm trace}\Big[\partial_{\xi_n}\pi^+_{\xi_n}\sigma_{-1}(\nabla^{\Psi}_{X}\nabla^{\Psi}_{Y}D_{\Psi}^{-1})\times
2\Big(\frac{c(\xi)c(\Psi)c(\xi)}{|\xi|^6}-\frac{c(\Psi)}{|\xi|^4}\Big)\Big]\nonumber\\
&=\sum_{j,l=1}^{n-1}X_jY_l\xi_j\xi_l\frac{-i\pi}{2}{\rm trace}[c(dx_n)c(\Psi))]
+X_nY_n\frac{i\pi}{2}{\rm trace}[c(dx_n)c(\Psi))].
\end{align}
From (\ref{61}),(\ref{71}) and (\ref{710}), we get
\begin{align}\label{74}
\widetilde{\Phi}_5
=&\sum_{j,l=1}^{n-1}X_jY_l\frac{-2\pi^{2}}{3}{\rm trace}[c(dx_n)c(\Psi))]\Omega_{3}dx'
+X_nY_n\frac{i\pi}{2}{\rm trace}[c(dx_n)c(\Psi))]\Omega_{3}dx'\nonumber\\
&+\sum_{j,l=1}^{n-1}X_jY_lh'(0)(\frac{55}{26}+\frac{85i}{24})\pi^2\Omega_{3}dx'+
X_nY_n(-\frac{50+7i}{16}) \pi\Omega_{3}dx'
\end{align}

Let $X=X^T+X_n\partial_n,~Y=Y^T+Y_n\partial_n,$ then we have $\sum\limits_{j=1}^{n-1}X_jY_j=g(X^T,Y^T).$
 Now $\widehat{\Phi}$ is the sum of the cases (a), (b) and (c). Combining with the five cases, this yields
\begin{align}\label{1795}
\widetilde{\Phi}=\sum_{i=1}^5\widetilde{\Phi}_i
=&\big[(\frac{-4681}{24}+\frac{5i}{6})\pi^2+(\frac{55}{26}
+\frac{85i}{24})\big]g(X^T,Y^T)h'(0)\Omega_3dx'+\big[(\frac{451}{4}+7i+\frac{\pi}{8}) h'(0)\nonumber\\
&-\frac{50+7i}{16} \big]X_nY_n\pi\Omega_3dx'-\big[\frac{3\pi+2\pi i}{4}X_nY_n+\pi^2g(X^T,Y^T) \big]{\rm trace}[c(dx_n)c(\Psi)]\Omega_3dx'\nonumber\\
&-\frac{9}{16}h'(0)\pi\Omega_3dx'+\frac{3\pi i}{16}{\rm trace}
[-c(X)c(\Psi)X_n-c(Y)c(\Psi)Y_n]\Omega_3dx'
.
\end{align}
So, we are reduced to prove the following.

Combining (\ref{fum4.3}) and (\ref{1795}), we obtain Theorem 1.2.\\

When $c(\Psi)=f$, we can calculate
${\rm trace}~[c(dx_n)f]=0$ and ${\rm trace}~[c(X)f]=0$,
then we can straightforwardly assert the following corollary:
\begin{cor}
Let $M$ be a $4$-dimensional oriented
compact manifolds with the boundary $\partial M$ and the metric
$g^M$ as above, and let $c(\Psi)=f,$ then
\begin{align}
&\widetilde{{\rm Wres}}[\pi^{+}(\nabla^{\Psi}_{X}\nabla^{\Psi}_{Y}D_{\Psi}^{-1})
    \circ\pi^{+}(D_{\Psi}^{-3})]=\frac{4\pi^2}{3}\int_{M}EG(X,Y)vol_{g}
+\frac{1}{2}\int_{M}{\rm trace}[\frac{1}{4}s-3f^2]g(X,Y)dvol_{M}\nonumber\\&
+\int_{\partial M}\bigg\{\big[(\frac{-4681}{24}+\frac{5i}{6})\pi^2+(\frac{55}{26}
+\frac{85i}{24})\big]g(X^T,Y^T)h'(0)\Omega_3+\big[(\frac{451}{4}+7i+\frac{\pi}{8}) h'(0)\nonumber\\
&-\frac{50+7i}{16} \big]X_nY_n\pi\Omega_3-\frac{9}{16}h'(0)\pi\Omega_3\bigg\}d{\rm Vol_{\partial M}},
\end{align}
where $s$ is the scalar curvature.
\end{cor}

By (\ref{4522}), when $c(\Psi)=c(\overline{X})$, 
then we can straightforwardly assert the following corollary:
\begin{cor}
Let $M$ be a $4$-dimensional oriented
compact manifolds with the boundary $\partial M$ and the metric
$g^M$ as above, and let $c(\Psi)=c(\overline{X}),$ then
\begin{align}
&\widetilde{{\rm Wres}}[\pi^{+}(\nabla^{\Psi}_{X}\nabla^{\Psi}_{Y}D_{\Psi}^{-1})
    \circ\pi^{+}(D_{\Psi}^{-3})]=\frac{4\pi^2}{3}\int_{M}EG(X,Y)vol_{g}
-4\pi^2\int_{M}\Big[g(X,\nabla^{TM}_Y\overline{X})
+g(Y,\nabla^{TM}_X\overline{X})\Big]dvol_{M}\nonumber\\
&+\frac{1}{2}\int_{M}{\rm trace}\big[\frac{1}{4}s+\sum\limits^{n}_{j=1}\frac{1}{2}c(\overline{X})c(e_j)
c(\overline{X})c(e_j)+|\overline{X}|^2\big]g(X,Y)dvol_{M}
+\int_{\partial M}\bigg\{\big[(\frac{-4681}{24}+\frac{5i}{6})\pi^2+(\frac{55}{26}
+\frac{85i}{24})\big]\nonumber\\
&\times g(X^T,Y^T)h'(0)\Omega_3
+\big[(\frac{451}{4}+7i+\frac{\pi}{8}) h'(0)-\frac{50+7i}{16} \big]X_nY_n\pi\Omega_3+\big[(3\pi+2\pi i)X_nY_n+\pi^2g(X^T,Y^T) \big]
\nonumber\\
&\times g(\partial x_n,X)\Omega_3-\frac{9}{16}h'(0)\pi\Omega_3dx'+\frac{3\pi i}{16}{\rm trace}
[4g(X,\overline{X})X_n+4g(Y,\overline{X})Y_n]\Omega_3
\bigg\}d{\rm Vol_{\partial M}},
\end{align}
where $s$ is the scalar curvature.
\end{cor}
By (\ref{4555}), when $c(\Psi)=c(X_1)c(X_2),$
we compute $\widetilde{{\rm Wres}}[\pi^{+}(\nabla^{\Psi}_{X}\nabla^{\Psi}_{Y}D_{\Psi}^{-1})
    \circ\pi^{+}(D_{\Psi}^{-3})].$
\begin{cor}
Let $M$ be a $4$-dimensional oriented
compact manifolds with the boundary $\partial M$ and the metric
$g^M$ as above, and let $c(\Psi)=c(X_1)c(X_2),$ then
\begin{align}
&\widetilde{{\rm Wres}}[\pi^{+}(\nabla^{\Psi}_{X}\nabla^{\Psi}_{Y}D_{\Psi}^{-1})
    \circ\pi^{+}(D_{\Psi}^{-3})]\nonumber\\
&=\frac{4\pi^2}{3}\int_{M}EG(X,Y)vol_{g}+\frac{1}{2}\int_{M}{\rm trace}[\frac{1}{4}s+\sum\limits^{n}_{j=1}\frac{1}{2}c(\Psi)c(e_j)
c(\Psi)c(e_j)-\big(c(\Psi)\big)^2]g(X,Y)dvol_{M}
\nonumber\\
&+\int_{\partial M}\bigg\{\big[(\frac{-4681}{24}+\frac{5i}{6})\pi^2+(\frac{55}{26}
+\frac{85i}{24})\big]g(X^T,Y^T)h'(0)\Omega_3+\big[(\frac{451}{4}+7i+\frac{\pi}{8}) h'(0)\nonumber\\
&-\frac{50+7i}{16} \big]X_nY_n\pi\Omega_3
-\frac{9}{16}h'(0)\pi\Omega_3dx'
\bigg\}d{\rm Vol_{\partial M}},
\end{align}
where $s$ is the scalar curvature.
\end{cor}

By (\ref{4588}), when $c(\Psi)=c(X_1)c(X_2)c(X_3),$ 
 we compute that:
\begin{cor}
Let $M$ be a $4$-dimensional oriented
compact manifolds with the boundary $\partial M$ and the metric
$g^M$ as above, and let $c(\Psi)=c(X_1)c(X_2)c(X_3),$ then
\begin{align}
&\widetilde{{\rm Wres}}[\pi^{+}(\nabla^{\Psi}_{X}\nabla^{\Psi}_{Y}D_{\Psi}^{-1})
    \circ\pi^{+}(D_{\Psi}^{-3})]
    =\frac{4\pi^2}{3}\int_{M}EG(X,Y)vol_{g}
+\frac{\pi^2}{2}\int_{M}\bigg\{
8\Big[
g(X,X_1)\Big(g(\nabla^{TM}_YX_2,X_3)\nonumber\\&+g(X_2,\nabla^{TM}_YX_3)\Big)
-g(X,X_2)\Big(g(\nabla^{TM}_YX_1,X_3)+g(X_1,\nabla^{TM}_YX_3)\Big)
+g(X,X_3)\Big(g(\nabla^{TM}_YX_1,X_2)\nonumber\\&+g(X_1,\nabla^{TM}_YX_2)\Big)
-g(X,\nabla^{TM}_YX_1)g(X_2,X_3)
-g(X,\nabla^{TM}_YX_2)g(X_1,X_3)
+g(X,\nabla^{TM}_YX_3)g(X_1,X_2)
\Big]\nonumber\\&
-
8\Big[
g(Y,X_1)\Big(g(\nabla^{TM}_XX_2,X_3)+g(X_2,\nabla^{TM}_XX_3)\Big)
-g(Y,X_2)\Big(g(\nabla^{TM}_XX_1,X_3)+g(X_1,\nabla^{TM}_XX_3)\Big)\nonumber\\&
+g(Y,X_3)\Big(g(\nabla^{TM}_XX_1,X_2)+g(X_1,\nabla^{TM}_XX_2)\Big)
-g(Y,\nabla^{TM}_YX_1)g(X_2,X_3)
-g(Y,\nabla^{TM}_YX_2)g(X_1,X_3)\nonumber\\&
+g(Y,\nabla^{TM}_YX_3)g(X_1,X_2)
\Big]dvol_{M}
+\int_{\partial M}\bigg\{\big[(\frac{-4681}{24}+\frac{5i}{6})\pi^2+(\frac{55}{26}
+\frac{85i}{24})\big]g(X^T,Y^T)h'(0)\Omega_3\nonumber\\
&+\big[(\frac{451}{4}+7i+\frac{\pi}{8}) h'(0)-\frac{50+7i}{16} \big]X_nY_n\pi\Omega_3-\big[\frac{3\pi+2\pi i}{4}X_nY_n+\pi^2g(X^T,Y^T) \big]\cdot\bigg(4g(\partial{x_n}, X_1)g(X_2, X_3)\nonumber\\
&-4g(\partial{x_n}, X_2)g(X_1, X_3)+4g(\partial{x_n}, X_3)g(X_1, X_2)\bigg)\Omega_3-\frac{9}{16}h'(0)\pi\Omega_3dx'
+\frac{3\pi i}{4}
\bigg[
\Big( g(X,X_1)g(X_2,X_3)\nonumber\\&-g(X,X_2)g(X_1,X_3)+g(X,X_3)g(X_1,X_2) \Big)X_n
+\Big( g(Y,X_1)g(X_2,X_3)-g(Y,X_2)g(X_1,X_3)
\nonumber\\&+g(Y,X_3)g(X_1,X_2) \Big)Y_n
\bigg]\Omega_3
\bigg\}d{\rm Vol_{\partial M}},
\end{align}
where $s$ is the scalar curvature.
\end{cor}

\section*{Acknowledgements}
This work was supported by NSFC No.12301063 and NSFC No.11771070, DUFE202159 and Basic research Project of the Education Department of Liaoning Province (Grant No. LJKQZ20222442). The authors thank the referee for his (or her) careful reading and helpful comments.

\end{document}